\documentclass[preprint,12pt]{elsarticle}

\usepackage{amssymb}
\usepackage{amsmath} 
\usepackage{graphicx}
\usepackage{hyperref}
\usepackage{booktabs}
\usepackage{listings}
\usepackage{xcolor}
\usepackage{ragged2e}
\usepackage{float}
\usepackage{algorithm} 
\usepackage{algpseudocode} 
\usepackage{caption}
\usepackage{setspace}
\usepackage[utf8]{inputenc}
\usepackage{lineno} 
\usepackage{xcolor}

\journal{}

\begin{document}

\begin{frontmatter}



\title{\textcolor{blue}{A polyhedral scaled boundary finite element method solving three-dimensional heat conduction problems}}


\author[inst1]{Mingjiao Yan}

\affiliation[inst1]{organization={College of Water Conservancy and Hydropower Engineering, Hohai University},
            city={Nanjing},
            postcode={210098}, 
            state={Jiangsu},
            country={China}}
\author[inst2]{Yang Yang}
\affiliation[inst2]{organization={PowerChina Kunming Engineering Corporation Limited},
            city={Kunming},
            postcode={650051}, 
            state={Yunnan},
            country={China}}
\affiliation[inst3]{organization={College of Water Conservancy, Yunnan Agricultural University},
            city={Kunming},
            postcode={650201}, 
            state={Yunnan},
            country={China}}
\author[inst1]{Chao Su}
\author[inst2]{Zongliang Zhang}
\author[inst3]{Qingsong Duan}
\author[inst1]{Dengmiao Hao}
\author[inst4]{Jian Zhou}
\affiliation[inst4]{organization={Water Conservancy Construction and Quality Safety Management Service Center of Zhenyuan Yi, Hani and Lahu Autonomous County},
            city={Pu'er},
            postcode={666599}, 
            state={Yunnan},
            country={China}}
\begin{abstract}
In this study, we derived a three-dimensional scaled boundary finite element formulation for heat conduction problems. By incorporating Wachspress shape functions, a polyhedral scaled boundary finite element method (PSBFEM) was proposed to address heat conduction challenges in complex geometries. To address the complexity of traditional methods, this work introduced polygonal discretization techniques that simplified the topological structure of the polyhedral mesh and effectively integrated polyhedral and octree meshes, thereby reducing the number of element faces and enhancing mesh efficiency to accommodate intricate shapes. The developed formulation supported both steady-state and transient heat conduction analyses and was implemented in ABAQUS through a user-defined element (UEL). Through a series of numerical examples, the accuracy and convergence of the proposed method were validated. The results indicated that the PSBFEM consistently achieved higher accuracy than the FEM as the mesh was refined. The polyhedral elements offered a computationally efficient solution for complex simulations, significantly reducing computational costs. \textcolor{blue}{Additionally, by utilizing the octree mesh parent element acceleration technique, the computational efficiency of PSBFEM surpassed that of the FEM.}
\end{abstract}

\begin{highlights}
\item Derived a semi-analytical formulation to solve heat conduction problems based on the SBFEM
\item Developed a polyhedral SBFEM (PSBFEM) framework incorporating Wachspress shape functions
\item Implemented the PSBFEM in ABAQUS UEL to solve both steady-state and transient heat conduction problems
\item Achieved superior accuracy with the PSBFEM compared to the FEM during mesh refinement
\item Polyhedral elements provide an efficient solution for complex simulations by significantly reducing computational costs
\end{highlights}


\end{frontmatter}


\section{Introduction}
\label{sec:Introduce}
Heat conduction problems are prevalent in modern engineering applications, such as electronic device cooling \cite{zhang2021review}, material processing \cite{murz2021}, nuclear reactor design \cite{tian2021experimental}, and civil engineering \cite{hunganhOver2021a}. Accurately and efficiently solving three-dimensional (3D) heat conduction problems is essential for optimizing designs and enhancing system performance. Experimentation is often employed to study heat conduction phenomena \cite{shiExp2010,yangControl2021}. However, these methods can be time-consuming and costly, potentially failing to capture all the nuances of the system, particularly in geometrically intricate or dynamically varying scenarios. Additionally, although analytical methods \cite{pop2020,hob2021} have been used over the past decades, they are typically limited by the complexities of geometry and boundary conditions. Consequently, numerical methods, particularly finite element method (FEM), have become indispensable for addressing complex heat conduction problems \cite{malek2020three,bilal2020finite}.

In traditional 3D FEM, tetrahedral and hexahedral elements are widely utilized. However, the accuracy and reliability of FEM solutions depend heavily on mesh quality. Tetrahedral elements allow for automated mesh generation but often result in lower accuracy, while hexahedral elements provide higher accuracy but require manual intervention during the mesh generation process, making it more time-consuming. To achieve precision in computational results, high-quality meshes are typically generated through complex and time-consuming re-meshing algorithms. According to Sandia National Laboratories, preprocessing steps such as mesh generation account for more than 80\% of the total analysis time \cite{cottrell2009isogeometric}. To alleviate the burden of preprocessing, researchers have explored various alternative methods, including the boundary element method (BEM) \cite{majchrzak2023bem}, isogeometric analysis (IGA) \cite{yu2020locally}, meshfree methods \cite{afrasiabi2020contemporary}, and physics-informed neural networks (PINNs) \cite{cai2021physics}. Each of these methods aims to address specific challenges posed by the FEM, particularly in terms of geometric complexity and computational cost, yet they present their own limitations in practice.

For large-scale problems, employing a locally refined mesh in regions of interest has become a common strategy, as it provides higher accuracy in critical areas while reducing the overall computational load \cite{li2020n}. This approach has given rise to non-matching mesh techniques, which allow for tailored meshing strategies that apply finer meshes where needed and coarser meshes elsewhere \cite{lacroix2024comparative, damirchi2021transverse}. However, non-matching meshes introduce the challenge of handling hanging nodes—nodes that do not align with those in adjacent elements. Traditional FEM struggles with such nodes due to the disruption of continuity and compatibility conditions, which are essential for accurate solutions \cite{ooi2015adaptation}.

Several techniques have been developed to manage hanging nodes, including the Arlequin method \cite{sun2017mixed}, the mortar segment-to-segment method \cite{zhou2020three}, multi-point constraint equations \cite{jendele2009solution}, and the polygonal finite element method (PFEM) \cite{biabanaki2012polygonal,chi2015polygonal}. Among these, the PFEM inherently addresses non-matching mesh problems by treating hanging nodes as part of the element discretization. However, the PFEM's reliance on a weak formulation in each discretization direction often results in reduced computational accuracy. 

The scaled boundary finite element method (SBFEM) has emerged as a robust technique for handling hanging nodes in quadtree and octree meshes \cite{gravenkamp2017efficient,yang2022novel}. The SBFEM combines the strengths of both the FEM and the BEM into a semi-analytical approach. Discretization occurs only in the circumferential direction, while the radial direction retains a strong form of the governing equations. This hybrid approach enables significant computational efficiency, making it well-suited for complex engineering applications \cite{song2018scaled}. The SBFEM's ability to handle irregular meshes and discontinuities while maintaining high precision has led to its widespread application in various fields \cite{song2018review,ye2021free,chen2018efficient}. 

Moreover, to enhance the capability of SBFEM in handling complex geometries, researchers have introduced polygonal and polyhedral techniques. Ye et al. \cite{ye2021free,yePSBFEM2021} integrated standard SBFEM with a polygonal mesh technique to address static and dynamic problems. \textcolor{blue}{Natarajan et al. \cite{natarajan2017scaled} developed an SBFEM for 3D convex polyhedra by employing surface discretization with Wachspress interpolants, ensuring high accuracy in linear elasticity. Furthermore, to accurately capture the complex geometries of polygonal and quadtree elements, researchers \cite{natarajan2017scaled,li2024improved,ZANG2023279} introduced a NURBS-based polygonal SBFEM by incorporating non-uniform rational B-splines (NURBS), enabling precise representation of curved boundaries. These advancements have significantly enhanced the flexibility and accuracy of SBFEM, establishing it as a powerful tool for various engineering applications.}

\textcolor{blue}{Currently, no dedicated software exists for polyhedral SBFEM. Researchers typically implement this method through secondary development based on commercial software. Ya et al. \cite{yaOpen2021} extended SBFEM to address interfacial problems using polyhedral meshes in ABAQUS UEL. Yang et al. \cite{yangDevelopment2020} developed UEL and VUEL implementations in ABAQUS to perform static and elasto-dynamic stress analyses of solids. However, these approaches discretize the boundary surfaces using quadrilateral and triangular elements, resulting in polyhedral elements with relatively complex topologies.}

Recent studies have extended SBFEM to heat conduction problems. For instance, Bazyar et al. \cite{bazyar2015scaled} used the SBFEM to solve 2D heat conduction problems in anisotropic media. Yu et al. \cite{yu2021scaled} applied a hybrid quadtree mesh in SBFEM to analyze transient heat conduction problems involving cracks or inclusions. Yang et al. \cite{yang2021polygonal} developed a user element (UEL) of polygonal SBFEM to solve 2D heat conduction problems. However, most current research focuses on 2D problems. \textcolor{blue}{He et al. \cite{he2019image} presented a novel numerical technique that combines the quadtree technique, SBFEM, image-based modeling, and inverse analysis to evaluate the effective thermal conductivity of heterogeneous materials. He et al. \cite{he2019new} proposed a new adaptive algorithm integrating quadtree SBFEM and the smoothed effective heat capacity method for phase-change heat transfer problems. Wang et al. \cite{wang2023multiscale} developed a numerical model using SBFEM-based basis functions to bridge small and large scales in heat conduction, ensuring both efficiency and accuracy for complex geometries.} For 3D heat conduction problems, Lu et al. \cite{lu2017modified} developed a modified SBFEM to address steady-state heat conduction problems, but transient analysis remains underexplored. 

In this study, a polyhedral SBFEM (PSBFEM) framework is proposed for solving 3D steady-state and transient heat conduction problems. The remainder of this paper is organized as follows: Sections \ref{sec:2} and \ref{sec:3} derive the formulation of the SBFEM for 3D heat conduction analysis. Section \ref{sec:4} develops a PSBFEM by incorporating Wachspress shape functions. Section \ref{sec:5} presents the solution procedure for the SBFEM equations. Section \ref{sec:6} describes the implementation process for solving heat conduction problems using the SBFEM in ABAQUS UEL. Section \ref{sec:7} provides several numerical examples to demonstrate the accuracy and efficiency of the proposed method. Finally, Section \ref{sec:8} concludes the study by summarizing the key findings.

\section{Governing equations of the 3D heat conduction problems}
\label{sec:2}
In this study, we considered a 3D transient heat conduction problem, the governing equations without heat sources are written as \cite{lu2017modified,songScaled1999}
\textcolor{blue}{
\begin{equation}
    \mathbf{L}^\mathrm{T}\mathbf{q}+ \rho c\dot{{T}}=0\quad\text{in } \quad\Omega,  \label{eq:gov} 
\end{equation}}
where $\rho$ denotes the mass density. $c$ is the heat capacity, $\dot{T}$ is the temperature change rate, $\Omega$ denotes the computational domain; $\mathbf{q}$ is the heat flux vector, $\mathbf{q}$ can be expressed as
\begin{equation}
    \mathbf{q}=-\mathbf{k}\mathbf{L}T,
\end{equation}
where $T$ is the temperature, $\mathbf{k}$ is the matrix of the thermal conductivity, $\mathbf{k}$ can be written as:
\begin{equation}\mathbf{k}=\begin{bmatrix}k_x&0&0\\0&k_y&0\\0&0&k_z\end{bmatrix},\end{equation}
\textcolor{blue}{where $k_x$, $k_y$, and $k_z$ are the thermal conductivity in the x, y, and z directions, respectively. The operator $\mathbf{L}$ is the differential operator and can be written as}
\textcolor{blue}{\begin{equation}
    \mathbf{L}=\begin{Bmatrix}\frac{\partial}{\partial\hat{x}}\\\frac{\partial}{\partial\hat{y}}\\\frac{\partial}{\partial\hat{z}}\end{Bmatrix}.
\end{equation}}

The initial conditions can be written as follows:
\begin{equation}T(x,y,z,t=0)=T_0(x,y,z)\quad\text{in }\quad\Omega. \end{equation}
\textcolor{blue}{where $T_0$ denotes the initial temperature.} And the boundary conditions can be expressed as:
\begin{equation}T(x,y,z)=\bar{T}\quad\mathrm{on~}S_1,\end{equation}
\begin{equation}-k\frac{\partial T}{\partial n}=\overline{q}_n=\begin{cases} \quad q_2&\text{on}\quad S_2\\h(T-T_\infty)&\text{on}\quad S_3\end{cases},\end{equation}
where $n$ is outward normal vector to $\Omega$. $S$ is the boundary. In addition, $\bar{T}$ and $q_2$ denote the specified boundary temperature and heat flux, respectively, while $h$ represents the convection heat transfer coefficient, and $T_\infty$ is the surrounding temperature.

Applying the Fourier transform to Eq. (\ref{eq:gov}) yields the governing equation in the frequency domain as follows:
\textcolor{blue}{\begin{equation}\mathbf{L}^\mathrm{T}\mathbf{\tilde{q}}+\mathrm{i}\omega \rho c \tilde{T} =0
, \label{eq:threegov}\end{equation}}
where $\mathbf{\tilde{q}}$ and $\tilde{T}$ are the Fourier transform of $\mathbf{q}$ and $T$, respectively. $\omega$ is the frequency. When $\omega = 0$, the problem reduces to a steady-state heat conduction scenario.
\section{Formulation of SBFEM for heat conduction problems}
\label{sec:3}
\subsection{Scaled boundary transformation of the geometry}
\begin{figure}[H]
  \centering
  \includegraphics[width=0.5\textwidth]{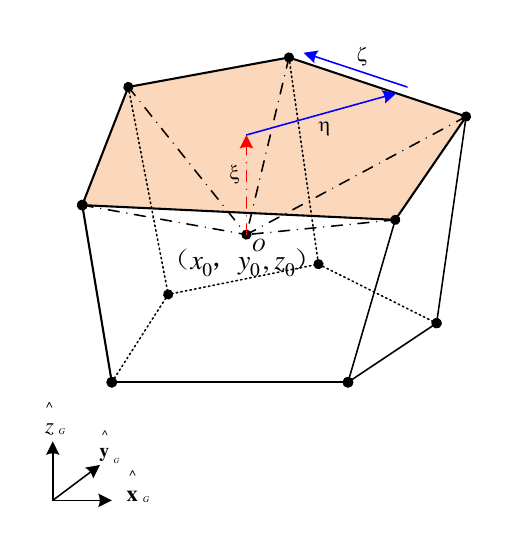}
  \caption{Scaled boundary coordinates ($\xi, \eta, \zeta$) representing the coordinate system centered at $O$ for an arbitrary faceted polyhedron.}
  \label{fig:SBFEM coordinates}
\end{figure}

The scaled boundary coordinates is shown in Fig. \ref{fig:SBFEM coordinates}.
The Cartesian coordinates of a point $(\hat{x},\hat{y},\hat{z})$ within the volume sector $\mathrm{V}_{\mathrm{e}}$ can be defined by the scaled boundary coordinates $(\xi,\eta,\zeta)$ as \cite{song2004matrix,song2018scaled}
\begin{subequations}
\begin{align}
    \hat{x}(\xi,\eta,\zeta)=\xi\hat{x}(\eta,\zeta)=\xi\mathbf{N}(\eta,\zeta)\mathbf{\hat{x}}, \\
    \hat{y}(\xi,\eta,\zeta)=\xi\hat{y}(\eta,\zeta)=\xi\mathbf{N}(\eta,\zeta)\mathbf{\hat{y}}, \\
    \hat z(\xi,\eta,\zeta)=\xi\hat z(\eta,\zeta)=\xi\mathbf{N}(\eta,\zeta)\mathbf{\hat z},
\end{align}
\end{subequations}
where $\mathbf{\hat{x}}$, $\mathbf{\hat{y}}$ and $\mathbf{\hat{z}}$ denote the nodal coordinate vectors of the surface element $\mathrm{S_e}$ in Cartesian coordinates. $\mathbf{N}(\eta,\zeta)$ is the shape function vector
\begin{equation}
    \mathbf{N}(\eta,\zeta)=[ N_1(\eta,\zeta) \quad N_2(\eta,\zeta) \quad ... \quad  N_n(\eta,\zeta)],
\end{equation}
where $N_n(\eta,\zeta)$ denotes the shape function and $n$ represents the total number of nodes in the polygon.

The partial derivatives with respect to the scaled boundary coordinates are connected to those in Cartesian coordinates by the following equation:
\begin{equation}\begin{Bmatrix}\frac{\partial}{\partial\xi}\\\frac{\partial}{\partial\eta}\\\frac{\partial}{\partial\zeta}\end{Bmatrix}=\begin{bmatrix}1&0&0\\0&\xi&0\\0&0&\xi\end{bmatrix}\mathbf{J_b}(\eta,\zeta)\begin{Bmatrix}\frac{\partial}{\partial \hat{x}}\\\frac{\partial}{\partial \hat{y}}\\\frac{\partial}{\partial \hat{z}}\end{Bmatrix}, \label{eq:L1}\end{equation}
where $\mathbf{J_b}(\eta,\zeta)$ is the Jacobian matrix on the boundary $(\xi=1)$, and is defined as follows: 
\begin{equation}\mathbf{J_b}(\eta,\zeta)=\begin{bmatrix}\hat{x}(\eta,\zeta)&\hat{y}(\eta,\zeta)&\hat{z}(\eta,\zeta)\\\hat{x}(\eta,\zeta)_{,\eta}&\hat{y}(\eta,\zeta)_{,\eta}&\hat{z}(\eta,\zeta)_{,\eta}\\\hat{x}(\eta,\zeta)_{,\zeta}&\hat{y}(\eta,\zeta)_{,\zeta}&\hat{z}(\eta,\zeta)_{,\zeta}\end{bmatrix}.\end{equation}

The determinant of $\mathbf{J_b}(\eta,\zeta)$ is given by the following expression:
\begin{equation}|\mathbf{J_b}|=\hat{x}(\hat{y}_{,\eta}\hat{z}_{,\zeta}-\hat{z}_{,\eta}\hat{y}_{,\zeta})+\hat{y}(\hat{z}_{,\eta}\hat{x}_{,\zeta}-\hat{x}_{,\eta}\hat{z}_{,\zeta})+\hat{z}(\hat{x}_{,\eta}\hat{y}_{,\zeta}-\hat{y}_{,\eta}\hat{x}_{,\zeta}),\end{equation}
where the argument $(\eta,\zeta)$ are omitted for simplicity. From Eq. (\ref{eq:L1}), it follows that the inverse relationship is expressed as:
\begin{equation}\left\{\begin{array}{c}\frac{\partial}{\partial\hat{x}}\\\frac{\partial}{\partial\hat{y}}\\\frac{\partial}{\partial\hat{z}}\end{array}\right\}=\mathbf{J}_{\mathbf{b}}^{-1}\left(\eta,\zeta\right)\left\{\begin{array}{c}\frac{\partial}{\partial\xi}\\\frac{1}{\xi}\frac{\partial}{\partial\eta}\\\frac{1}{\xi}\frac{\partial}{\partial\zeta}\end{array}\right\},\label{eq:L2}\end{equation}
where the inverse of the Jacobian matrix at the boundary, $\mathbf{J_b}^{-1}(\eta, \zeta)$, can be expressed as follows:
\begin{equation}\left.\mathbf{J}_{\mathbf{b}}^{-1}\left(\eta,\zeta\right)=\frac1{\left|\mathbf{J_b}\right|}\left[\begin{array}{ccc}y_{,\eta}z_{,\zeta}-z_{,\eta}y_{,\zeta}&zy_{,\zeta}-yz_{,\zeta}&yz_{,\eta}-zy_{,\eta}\\\\z_{,\eta}x_{,\zeta}-x_{,\eta}z_{,\zeta}&xz_{,\zeta}-zx_{,\zeta}&zx_{,\eta}-xz_{,\eta}\\\\x_{,\eta}y_{,\zeta}-y_{,\eta}x_{,\zeta}&yx_{,\zeta}-xy_{,\zeta}&xy_{,\eta}-yx_{,\eta}\end{array}\right.\right].\end{equation}

Hence, Eq. (\ref{eq:L2}) can also be expressed as
\begin{equation}\begin{aligned}\left\{\begin{array}{c}\frac{\partial}{\partial\hat{x}}\\\frac{\partial}{\partial\hat{y}}\\\frac{\partial}{\partial\hat{z}}\end{array}\right\}&=\frac{1}{|\mathbf{J}_{\mathbf{b}}|}\left[\begin{array}{c}y_{,\eta}z_{,\zeta}-z_{,\eta}y_{,\zeta}\\\\z_{,\eta}x_{,\zeta}-x_{,\eta}z_{,\zeta}\\\\x_{,\eta}y_{,\zeta}-y_{,\eta}x_{,\zeta}\end{array}\right]\frac{\partial}{\partial\xi}\\&+\frac{1}{\xi}\left(\frac{1}{|\mathbf{J}_{\mathbf{b}}|}\left[\begin{array}{c}zy_{,\zeta}-yz_{,\zeta}\\\\xz_{,\zeta}-zx_{,\zeta}\\\\yx_{,\zeta}-xy_{,\zeta}\end{array}\right]\frac{\partial}{\partial\eta}+\frac{1}{|\mathbf{J}_{\mathbf{b}}|}\left[\begin{array}{c}yz_{,\eta}-zy_{,\eta}\\\\zx_{,\eta}-xz_{,\eta}\\\\xy_{,\eta}-yx_{,\eta}\end{array}\right]\frac{\partial}{\partial\zeta}\right).\end{aligned}\end{equation}

The gradient operator in the Cartesian coordinate system can be converted to the scaled boundary coordinate system as follows:
\begin{equation}\mathbf{L}=\mathbf{b}_{1}\left(\eta,\zeta\right)\frac{\partial}{\partial\xi}+\frac{1}{\xi}\left(\mathbf{b}_{2}\left(\eta,\zeta\right)\frac{\partial}{\partial \eta}+\mathbf{b}_{3}\left(\eta,\zeta\right)\frac{\partial}{\partial\zeta}\right),\label{eq:L3}\end{equation}
where $\mathbf{b}_1\left(\eta,\zeta\right)$, $\mathbf{b}_2\left(\eta,\zeta\right)$, and $\mathbf{b}_3\left(\eta,\zeta\right)$ can be defined as 
\begin{equation}
\mathbf{b}_1\left(\eta,\zeta\right)=\left.\frac1{|\mathbf{J_b}|}\left[\begin{array}{c}y_{,\eta} z_{,\zeta}-z_{,\eta} y_{,\zeta}\\\\z_{,\eta} x_{,\zeta}-x_{,\eta} z_{,\zeta}\\\\x_{,\eta} y_{,\zeta}-y_{,\eta} x_{,\zeta}\end{array}\right.\right],
\end{equation}
\begin{equation}\mathbf{b}_2\left(\eta,\zeta\right)=\left.\frac{1}{|\mathbf{J_b}|}\left[\begin{array}{c}zy_{,\zeta}-yz_{,\zeta}\\\\xz_{,\zeta}-zx_{,\zeta}\\\\yx_{,\zeta}-xy_{,\zeta}\end{array}\right.\right],\end{equation}
\begin{equation}\mathbf{b}_3\left(\eta,\zeta\right)=\left.\frac{1}{|\mathbf{J_b}|}\left[\begin{array}{c}yz_{,\eta}-zy_{,\eta}\\\\zx_{,\eta}-xz_{,\eta}\\\\xy_{,\eta}-yx_{,\eta}\end{array}\right.\right].\end{equation}

\subsection{Temperature field}
The temperature of any point  $\tilde{T}(\xi,\eta,\zeta)$ in the SBFEM coordinates can be written as:
\begin{equation}
    \tilde{T}(\xi,\eta,\zeta)=\mathbf{N}(\eta,\zeta)\tilde{T}(\xi), \label{eq:L4}
\end{equation}
where $\tilde{T}(\xi)$ is the nodal temperature vector and $\mathbf{N}(\eta,\zeta)$ is the matrix of shape function.

By using Eqs. (\ref{eq:L4}) and (\ref{eq:L3}), yields the temperature field partial derivatives as follows:
\textcolor{blue}{\begin{equation}\begin{bmatrix}
    \frac{\partial \tilde{T}}{\partial\hat{x}}&\frac{\partial \tilde{T}}{\partial \hat{y}}&\frac{\partial \tilde{T}}{\partial \hat{z}}
\end{bmatrix}^\mathrm{T}=\mathbf{b}_{1}\left(\eta,\zeta\right)\frac{\partial \tilde{T}}{\partial\xi}+\frac{1}{\xi}\left(\mathbf{b}_{2}\left(\eta,\zeta\right)\frac{\partial \tilde{T}}{\partial \boldsymbol{\eta}}+\mathbf{b}_{3}\left(\eta,\zeta\right)\frac{\partial \tilde{T}}{\partial\boldsymbol{\zeta}}\right).\label{eq:L5}\end{equation}}

For simplicity, Eq. (\ref{eq:L5}) can be rewritten by substituting the temperature expression from Eq. (\ref{eq:L4}), as follows:
\begin{equation}
\begin{bmatrix}\frac{\partial \tilde{T}}{\partial\hat{x}}&\frac{\partial \tilde{T}}{\partial\hat{y}}&\frac{\partial \tilde{T}}{\partial\hat{z}}\end{bmatrix}^\mathrm{T}=\mathbf{B}_1\left(\eta,\zeta\right)\tilde{T}_{,\xi}+\frac1\xi \mathbf{B}_2\left(\eta,\zeta\right) \tilde{T},
\end{equation}
where 
\begin{equation}
\mathbf{B}_1\left(\eta,\zeta\right)=\mathbf{b}_{1}\left(\eta,\zeta\right)\mathbf{N}(\eta,\zeta),
\end{equation}
\begin{equation}
\mathbf{B}_2\left(\eta,\zeta\right)=\mathbf{b}_{2}\left(\eta,\zeta\right)\mathbf{N}(\eta,\zeta)_{,\eta}+\mathbf{b}_{3}\left(\eta,\zeta\right)\mathbf{N}(\eta,\zeta)_{,\zeta}.
\end{equation}

The heat flux $\mathbf{\tilde{q}}(\xi,\eta,\zeta)$ can be represented in the coordinate system of the SBFEM as follows:
\begin{equation}\mathbf{\tilde{q}}(\xi,\eta,\zeta)=-\mathbf{k}\Big(\mathbf{B_1}(\eta,\zeta) \tilde{T}(\xi)_{,\xi}+\frac{1}{\xi}\mathbf{B_2}(\eta,\zeta) \tilde{T}(\xi) \Big).\label{eq:heat flux}\end{equation}

\subsection{Scaled boundary finite element equation}
\textcolor{blue}{By employing the method of weighted residuals \cite{songScaled1999}, Eq. (\ref{eq:threegov}) can be transformed into} 
\textcolor{blue}{
\begin{equation}
\begin{aligned}&\int_\Omega w\mathbf{b_1}^\mathrm{T}\mathbf{\tilde{q},_\xi}\mathrm{~d}\Omega+\int_\Omega w\frac{1}{\xi}(\mathbf{b_2}^\mathrm{T}\mathbf{\tilde{q},_\eta}+\mathbf{b_3}^\mathrm{T}\mathbf{\tilde{q},_\zeta})\mathrm{~d}\Omega\\&
+\mathrm{i}\omega\int_\Omega w\rho c\tilde{T}\mathrm{d}\Omega=0.\end{aligned}\label{eq:weighted gov}
\end{equation}}
\textcolor{blue}{where $w = w(\xi, \eta, \zeta)$ is the weighting function. Following the procedures outlined by \cite{songScaled1999}, we further simplify Eq. (\ref{eq:weighted gov}) as follows (see \ref{appendix:A} for the detailed derivation):}
\textcolor{blue}{
\begin{align}
    &\mathbf{E}_0 \xi^2 \tilde{T}(\xi)_{,\xi\xi} 
    + \left( 2\mathbf{E}_0 - \mathbf{E}_1 + \mathbf{E}_1^{T} \right) \xi \tilde{T}(\xi)_{,\xi} \notag \\
    &+ \left( \mathbf{E}_1^\mathrm{T} - \mathbf{E}_2 \right) \tilde{T}(\xi) 
    - \mathbf{M}_0 \mathrm{i}\omega \xi^2 \tilde{T}(\xi) = \xi \mathbf{F}(\xi),  
    \label{eq:mainequation1}
\end{align}}
where the global coefficient matrices $\mathbf{E_0}$, $\mathbf{E_1}$, $\mathbf{E_2}$ and $\mathbf{M_0}$ for the entire element are formed by assembling the local coefficient matrices $\mathbf{E_{0}^{e}}$, $\mathbf{E_{1}^{e}}$, $\mathbf{E_{2}^{e}}$ and $\mathbf{M_{0}^{e}}$ corresponding to each surface element. These local matrices for a surface element $\mathrm{S_e}$ can be expressed as follows:
\begin{equation}\mathbf{E}_0^\mathrm{e}=\int_{\mathbf{S}_\mathrm{e}}\mathbf{B}_1^\mathrm{T}\mathbf{kB}_1\left|\mathbf{J}_\mathrm{b}\right|\mathrm{d}\eta\mathrm{d}\zeta,\label{eq:E0}\end{equation}
\begin{equation}\mathbf{E}_1^\mathrm{e}=\int_{\mathbf{S}_\mathrm{e}}\mathbf{B}_2^\mathrm{T}\mathbf{kB}_1\left|\mathbf{J}_\mathrm{b}\right|\mathrm{d}\eta\mathrm{d}\zeta,\end{equation}
\begin{equation}\mathbf{E}_2^\mathrm{e}=\int_{\mathrm{S}_\mathrm{e}}\mathbf{B}_2^\mathrm{T}\mathbf{kB}_2\left|\mathbf{J}_\mathrm{b}\right|\mathrm{d}\eta\mathrm{d}\zeta,\end{equation}
\begin{equation}\mathbf{M_0^e}=\int_{\mathrm{S_e}}\mathbf{N^T}\rho c \mathbf{N}\left|\mathbf{J_b}\right|\mathrm{d}\eta\mathrm{d}\zeta.\label{eq:M0}\end{equation}

\section{Polyhedral SBFEM element}
\label{sec:4}
\subsection{Polyhedral element construction}
\begin{figure}[H]
  \centering
  \includegraphics[width=1.0\textwidth]{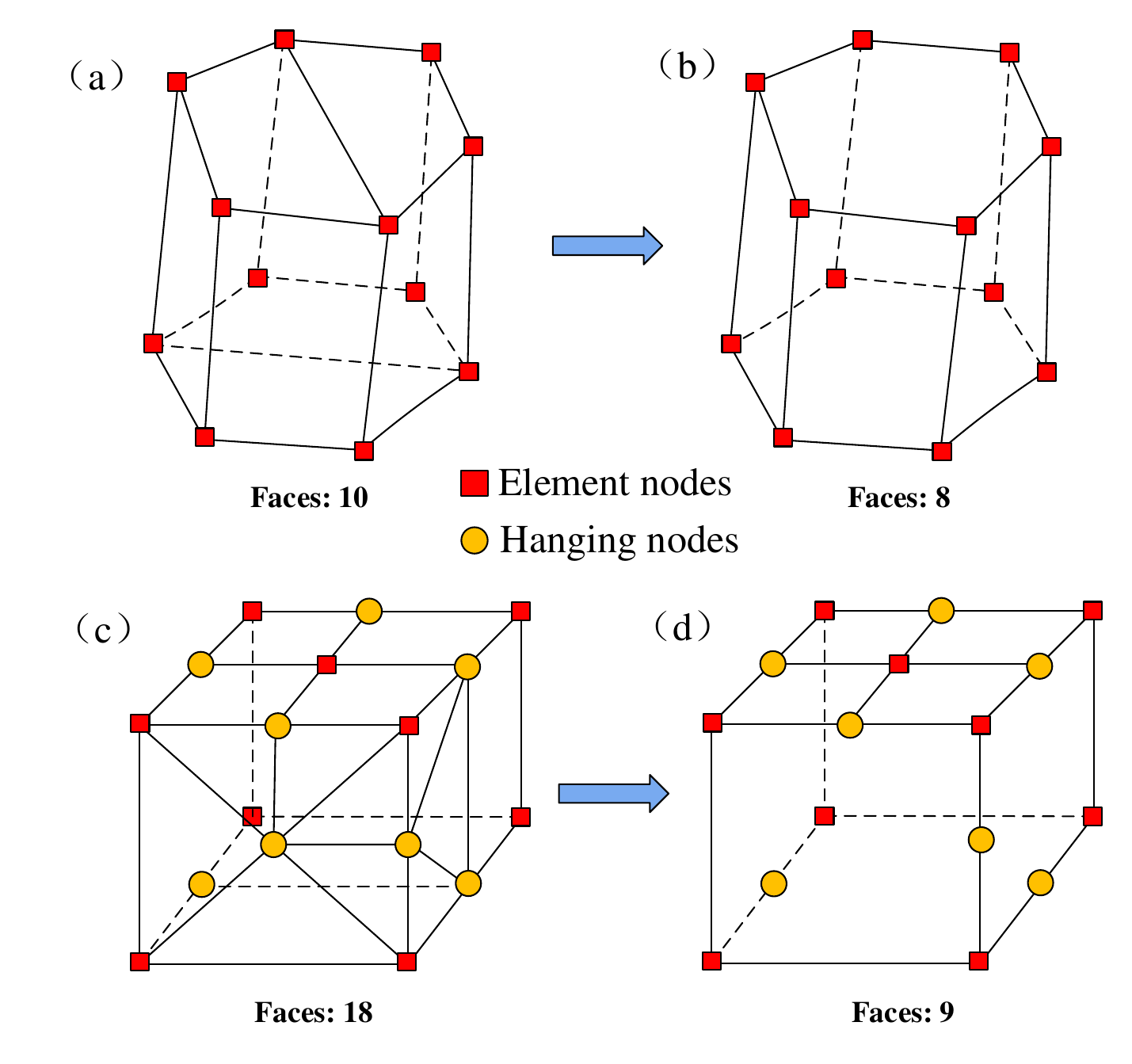}
  \caption{Polyhedral elements construction based on polygonal surfaces: (a) traditional polyhedral element; (b) polyhedral element constructed from polygonal surfaces; (c) traditional octree element; (d) octree element constructed from polygonal surfaces.}
  \label{fig:element construction}
\end{figure}

\begin{figure}[H]
  \centering
  \includegraphics[width=1.0\textwidth]{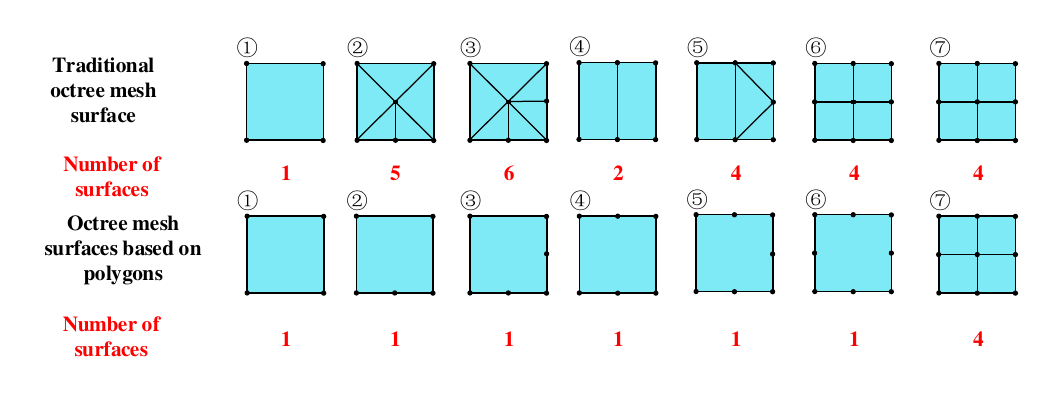}
  \caption{Comparison of the surfaces between traditional octree mesh and polygon-based octree mesh.}
  \label{fig:octree surface}
\end{figure}

In traditional 3D SBFEM, boundary surfaces are discretized using triangular and quadrilateral elements \cite{yaOpen2021,yangDevelopment2020}. As a result, the polyhedra processed by traditional 3D SBFEM exhibit relatively complex topologies, as illustrated in Figs. \ref{fig:element construction} (a) and (c). The intricate topological geometry of these polyhedra is evident. To address this complexity, this work introduces polygonal discretization techniques that simplify the topological structure of the polyhedra, thereby reducing the number of element faces, as shown in Figs. \ref{fig:element construction} (b) and (d). From the figures, it is clear that the polyhedra constructed using polygonal discretization significantly reduce the number of element faces, particularly for the octree mesh. This reduction not only enhances mesh efficiency but also provides clearer visualization. Furthermore, Fig. \ref{fig:octree surface} compares the traditional octree mesh with the polygon-based octree mesh surface, highlighting that the polygon-based octree mesh surface is more concise and has fewer faces.

\begin{figure}[H]
  \centering
  \includegraphics[width=0.8\textwidth]{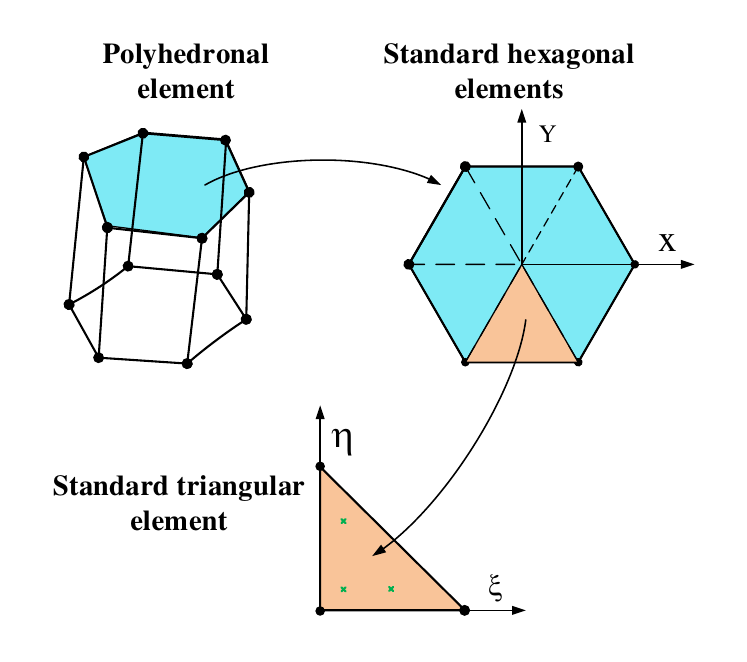}
  \caption{Numerical integration techniques for arbitrary polytopes.}
  \label{fig:Numerical integration technique}
\end{figure}

To construct polyhedral elements, we have introduced polygonal techniques. As shown in Fig. \ref{fig:Numerical integration technique}, the polygonal surface of the polyhedron is first mapped onto a regular polygon. This regular polygon is then subdivided into smaller triangles, and each of these sub-triangles is subsequently mapped onto a standard triangle. Quadrature rules, commonly used in the FEM for numerical integration, are applied to these sub-triangles to perform the integration. 

\subsection{Barycentric coordinates over polygonal surfaces}
\textcolor{blue}{
Various methods have been developed for constructing shape functions in polyhedral elements \cite{perumal2018brief}, and their comparisons can be found in \cite{floater2006general,sukumar2006recent}. The Wachspress shape function, which is based on a rational wedge function and derived using perspective geometry theory for collocation on convex polygon boundaries, is particularly notable for its applicability to arbitrary convex polygons \cite{wachspress2006rational}. The Wachspress shape function has several desirable properties, such as the partition of unity, Kronecker-delta function, non-negative, linear precision, etc \cite{nguyen2017polytree}. Moreover, the key advantage of the Wachspress shape function is its $C^\infty$ smoothness within the polygonal domain and $C^0$ continuity along the boundaries \cite{wu2023polygonal}. Hence, this property enhances the accuracy and ensures the convergence of temperature solutions.} 

\begin{figure}[H]
  \centering
  \includegraphics[width=0.65\textwidth]{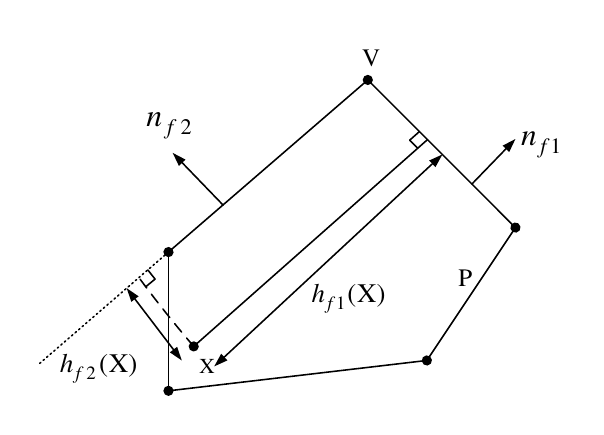}
  \caption{Barycentric coordinates: Wachspress basis function.}
  \label{fig:A ployhedral element}
\end{figure}

\textcolor{blue}{In this work, we focus on using Wachspress shape function to construct conforming approximations on polygonal meshing structures.} Wachspress introduced rational basis functions for polygonal elements based on principles from projective geometry. These functions ensure accurate nodal interpolation and maintain linearity along the boundaries by leveraging the algebraic equations of the edges. The polygonal element in Barycentric coordinates is presented in Fig. \ref{fig:A ployhedral element}. Let
\begin{equation}
w_k(\mathbf{x})=\frac{\operatorname{det}\left(\mathbf{n}_{f_1}, \mathbf{n}_{f_2}\right)}{h_{f_1}(\mathbf{x}) h_{f_2}(\mathbf{x})},
\end{equation}
where $\mathbf{n}_{f_1}$ and $\mathbf{n}_{f_2}$ denote the normal vectors of the boundary edges. $h_{f_1}(\mathbf{x})$ and $h_{f_2}(\mathbf{x})$ denote the perpendicular distances from $\mathbf{x}_k$ to the facets $f_1$ and $f_2$, which are calculated as follows:
\begin{equation}
h_{f_i}(\mathbf{x}) = (\mathbf{x}_k - \mathbf{x}) \cdot \mathbf{n}_{f_i}.
\end{equation}

The edges $f_i$ should be arranged in a counterclockwise order around the node $\mathbf{x}_k$, as viewed from the exterior. The shape functions are then given by \cite{warren2003uniqueness,warren2007barycentric}
\begin{equation}\phi_\mathbf{k}(\mathbf{x})=\frac{w_\mathbf{k}(\mathbf{x})}{\sum_{\mathbf{k}\in V}w_\mathbf{k}(\mathbf{x})}.\end{equation}

\section{Solution procedure}
\label{sec:5}
\subsection{Stiffness matrix solution}
\textcolor{blue}{In this study, the $\mathbf{F({\xi})}$ on the right side of Eq. (\ref{eq:mainequation1}) represents the influence of the normal flow rate across the side face. When the domain is enclosed or the side is insulated, $\mathbf{F({\xi})}$ becomes zero \cite{yu2021scaled,LI2016678}. Since we assume an adiabatic side or a closed domain in our analysis, we set $\mathbf{F({\xi})} = 0$.} The steady-state temperature field within the SBFEM formulation can be obtained by substituting $\omega=0$ into Eq. (\ref{eq:mainequation1}), yielding the following result:
\begin{equation}\mathbf{E}_0\xi^2\tilde{T}(\xi)_{,\xi\xi}+\left(2\mathbf{E}_0-\mathbf{E}_1+\mathbf{E}_1^ {T}\right)\xi\tilde{T}(\xi)_{,\xi}+\left(\mathbf{E}_1^\mathrm{T}-\mathbf{E}_2\right)\tilde{T}(\xi)=\mathbf{0}. \label{eq:mainequation}\end{equation}

By introducing the variable, 
\begin{equation}\mathbf{X}(\xi)=\begin{bmatrix}\xi^{0.5}\tilde{\mathbf{h}}(\xi)\\\xi^{-0.5}\tilde{\mathbf{q}}(\xi)\end{bmatrix},\end{equation}
the SBFEM equation can be reformulated as a first-order ordinary differential equation (ODE):
\begin{equation}\xi\mathbf{X}(\xi)_{,\xi}=\mathbf{Z}_\mathrm{p}\mathbf{X}(\xi),\end{equation}
where the coefficient matrix $\mathbf{Z}_\mathrm{p}$ is a Hamiltonian matrix. $\mathbf{Z}_\mathrm{p}$ is expressed as
\begin{equation}\mathbf{Z_p}=\begin{bmatrix}-\mathbf{E}_0^{-1}\mathbf{E}_1^\mathrm{T}+0.5\mathbf{I}&\mathbf{E}_0^{-1}\\\mathbf{E}_2-\mathbf{E}_1\mathbf{E}_0^{-1}\mathbf{E}_1^\mathrm{T}&\mathbf{E}_1\mathbf{E}_0^{-1}-0.5\mathbf{I}\end{bmatrix}.\label{eq:zp}\end{equation}

The eigenvalue decomposition of  $\mathbf{Z}_\mathrm{p}$ can be written as 
\begin{equation}\mathbf{Z}_\mathrm{p}\begin{bmatrix}\mathbf{\Phi}_\mathrm{h1}&&\mathbf{\Phi}_\mathrm{h2}\\\\\mathbf{\Phi}_\mathrm{q1}&&\mathbf{\Phi}_\mathrm{q2}\end{bmatrix}=\begin{bmatrix}\mathbf{\Phi}_\mathrm{h1}&&\mathbf{\Phi}_\mathrm{h2}\\\\\mathbf{\Phi}_\mathrm{q1}&&\mathbf{\Phi}_\mathrm{q2}\end{bmatrix}\begin{bmatrix}\boldsymbol{\Lambda}^+&&0\\\\0&&\boldsymbol{\Lambda}^-\end{bmatrix},\label{eq:eigen decomp}\end{equation}
where $\mathbf{\Phi}_{\mathrm{h1}}$ and $\mathbf{\Phi}_{\mathrm{q1}}$ are the corresponding eigenvectors of $\boldsymbol{\Lambda}^+$, and $\mathbf{\Phi}_{\mathrm{h2}}$ and $\mathbf{\Phi}_{\mathrm{q2}}$ are the corresponding eigenvectors of $\boldsymbol{\Lambda}^-$. $\boldsymbol{\Lambda}^+$ and $\boldsymbol{\Lambda}^+$ denote the diagonal matrices of eigenvalues with positive and negative real parts, respectively. As a whole, $\mathbf{\Phi}_{\mathrm{h}}$ and $\mathbf{\Phi}_{\mathrm{q}}$ is the modal temperature and flux, respectively. \textcolor{blue}{The solution of $\mathbf{X}(\xi)$ can be written as}

\begin{equation}\mathbf{X}(\xi)=\begin{bmatrix}\mathbf{\Phi}_\mathrm{h1}&&\mathbf{\Phi}_\mathrm{h2}\\\\\mathbf{\Phi}_{\mathrm{q1}}&&\mathbf{\Phi}_\mathrm{q2}\end{bmatrix}\begin{bmatrix}\xi^{\boldsymbol{\Lambda}^+}&&\mathbf{0}\\\\\mathbf{0}&&\xi^{\boldsymbol{\Lambda}^-}\end{bmatrix}\begin{Bmatrix}\mathbf{c}^{\mathrm{\mathrm{n}}}\\\\\mathbf{c}^{\mathrm{p}}\end{Bmatrix},\end{equation}
where $\mathbf{c}^{\mathrm{n}}$ and $\mathbf{c}^{\mathrm{p}}$ are the integration constants related to $\boldsymbol{\Lambda}^+$ and $\boldsymbol{\Lambda}^-$, respectively, and are determined from the boundary conditions. The solution for $\mathbf{\tilde{h}}(\xi)$ and $\mathbf{\tilde{q}}(\xi)$ can be expressed as
\begin{equation}\mathbf{\tilde{h}}(\xi)=\mathbf{\Phi}_{\mathrm{h1}}\xi^{-\boldsymbol{\Lambda}^+-0.5\mathbf{I}}\mathbf{c}^{\mathrm{n}},\end{equation}
\begin{equation}\mathbf{\tilde{q}(\xi)=\Phi}_{\mathrm{q}1}\xi^{-\boldsymbol{\Lambda}^++0.5\mathbf{I}}\mathbf{c}^{\mathrm{n}}.\end{equation}

The relationship between the nodal temperature functions and the nodal flux functions is derived by eliminating the integration constants, as shown below:
\begin{equation}\mathbf{\tilde{q}}(\xi)=\mathbf{\Phi}_{\mathrm{q1}}\mathbf{\Phi}_{\mathrm{h1}}^{-1}\xi\mathbf{\tilde{h}}(\xi).\end{equation}

On the boundary $(\xi=1)$, the nodal flux vector is defined as $\mathbf{Q} = \mathbf{\tilde{q}}(\xi=1)$, and the nodal temperature vector is given by $\mathbf{T} = \mathbf{\tilde{h}} (\xi=1)$. The relationship between the heat flux and temperature vectors is expressed as $\mathbf{Q} = \mathbf{K}\mathbf{T}$. Consequently, the stiffness matrix of the subdomain can be formulated as:
\begin{equation}\mathbf{K}=\mathbf{\Phi}_{\mathrm{q1}}\mathbf{\Phi}_{\mathrm{h1}}^{-1}.\label{eq:K}\end{equation}

\subsection{Mass matrix solution}
The mass matrix for a volume element is expressed as \cite{yang2022novel,song2018scaled}
\begin{equation}\mathbf{M}=\mathbf{\Phi}_{\mathrm{h}1}^{-\mathrm{T}}\int_0^1\xi^{\boldsymbol{\Lambda}^+}\mathbf{m}_0\xi^{\boldsymbol{\Lambda}^+}\xi\mathrm{d}\xi\boldsymbol{\Phi}_{\mathrm{h}1}^{-1},\label{eq:Massequation}\end{equation}
where the coefficient matrix $\mathbf{m}_0$ is defined as
\begin{equation}\mathbf{m}_0=\mathbf{\Phi}_{\mathrm{h}1}^\mathrm{T}\mathbf{M}_0\mathbf{\Phi}_{h1}.\end{equation}

By converting the integration component in Eq. (\ref{eq:Massequation}) into matrix form, the expression for the mass matrix can be rewritten as follows:

\begin{equation}\mathbf{M}=\mathbf{\Phi}_{\mathrm{h}1}^{-\mathrm{T}}\mathbf{m}\mathbf{\Phi}_{\mathrm{h}1}^{-1},\label{eq:M}\end{equation}
where 
\begin{equation}\mathbf{m}=\int_0^1\xi^{\boldsymbol{\Lambda}^+}\mathbf{m}_0\xi^{\boldsymbol{\Lambda}^+}\xi\mathrm{d}\xi.\end{equation}

Each element of the matrix $\mathbf{m}$ can be computed analytically, resulting in the following expression:
\begin{equation}m_{ij}=\frac{m_{0ij}}{\lambda_{ii}^++\lambda_{jj}^++2},\end{equation}
where $m_{0ij}$ denotes an entry of $\mathbf{m}_0$, and $\lambda^{+}_{ij}$ and $\lambda^{+}_{ii}$ represent specific entries of the matrix $\boldsymbol{\Lambda}^+$.

\subsection{Transient solution}
The nodal temperature relationship for a bounded domain can be formulated as a standard time-domain equation, utilizing the steady-state stiffness and mass matrices, as follows: 
\begin{equation} \mathbf{K}\mathbf{T}(t) + \mathbf{M}\dot{\mathbf{T}}(t) = \mathbf{Q}(t), \label{eq:L6} \end{equation} 
where the nodal temperature $\mathbf{T}(t)$ is the time-derivative continuous function. Solving for the functional solution in the time domain presents significant challenges. In this study, the backward difference method \cite{zienkiewicz1989finite} is employed to solve Eq. (\ref{eq:L6}). The time domain is discretized into discrete time steps, and the solution at each time node is incrementally computed based on the initial conditions. Subsequently, the nodal temperature at any given time is determined through interpolation.

At time $\left[t,t+\Delta t\right]$, the temperature change rate $\dot{\mathbf{T}}(t)$ can be expressed as
\begin{equation}\dot{\mathbf{T}}(t)=\frac{\Delta \mathbf{T}}{\Delta t}=\frac{\mathbf{T}(t)^{t+\Delta t}-\mathbf{T}(t)^t}{\Delta t}. \label{eq:L7}\end{equation}

Eq. (\ref{eq:L7}) is substituted into Eq. (\ref{eq:L6}), and the equation at time step $ t +\Delta t$ can be obtained as follows:
\begin{equation}\mathbf{K}^{t+\Delta t}+\frac{\mathbf{M}^{t+\Delta t}}{\Delta t}\mathbf{T}(t)^{t+\Delta t}=\mathbf{Q}(t)^{t+\Delta t}+\frac{\mathbf{M}^{t+\Delta t}}{\Delta t}\mathbf{T}(t)^t.\end{equation}

\section{Implementation}
\label{sec:6}
\subsection{\textcolor{blue}{Overview framework}}
\begin{figure}[H]
  \centering
  \includegraphics[width=1.0\textwidth]{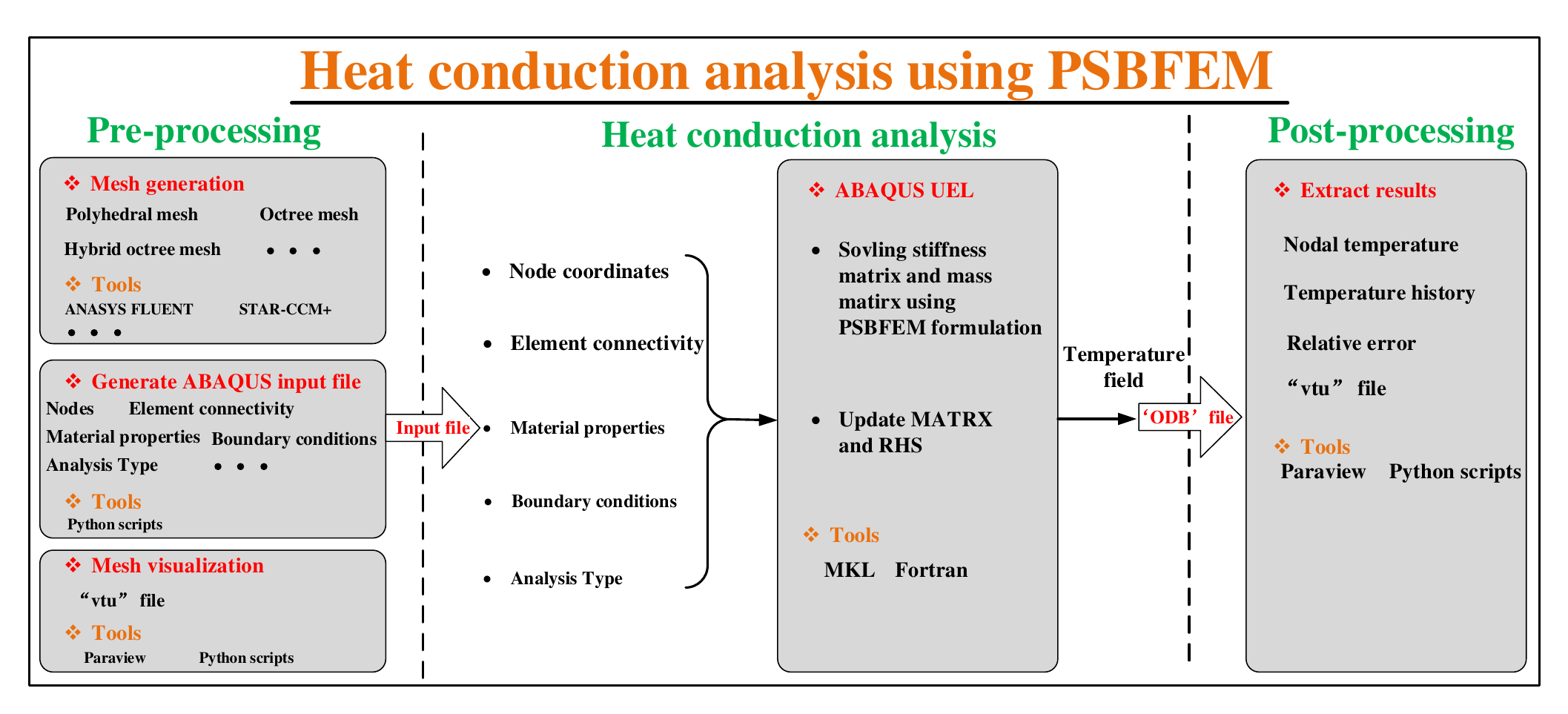}
  \caption{\textcolor{blue}{The framework of PSBFEM solving 3D heat conduction.}}
  \label{fig:flow}
\end{figure}

\textcolor{blue}{This study presents a framework for 3D heat conduction analysis using the PSBFEM with polyhedral and octree meshes, as illustrated in Fig. \ref{fig:flow}. The framework comprises three submodules: pre-processing, heat conduction analysis, and post-processing. In the pre-processing stage, commercial software such as ANSYS Fluent \cite{ANSYSFluent} and STAR-CCM+ \cite{STARCCM} is first used to generate polyhedral and octree meshes. A Python script is then employed to convert the mesh data into an ABAQUS input file (.inp). The 3D heat conduction analysis is conducted using an ABAQUS UEL subroutine implemented in FORTRAN. Further implementation details are provided in Section \ref{subsec:UEL}. Since ABAQUS CAE does not support the visualization of UEL elements, the results are extracted from the output database (.odb) file and visualized using ParaView \cite{ParaViewGuide}.}

\subsection{Implementation of the UEL}
\label{subsec:UEL}
Algorithm \ref{alg:1} provides a flowchart outlining the procedure for steady-state and transient heat conduction analysis, which has been successfully implemented in ABAQUS using the User Element (UEL) programming interface. The primary function of the UEL in ABAQUS is to update the element's contribution to the residual force vector (RHS) and the stiffness matrix (AMATRX) via the user subroutine interface provided by the software. For steady-state heat conduction analysis, the AMATRX and the RHS are defined as follows:
\begin{equation}
\mathrm{AMATRX}=\mathbf{K},\label{eq:amatrx-steady-state}
\end{equation}
\begin{equation}
    \mathrm{RHS}=-\mathbf{K}\mathbf{U},\label{eq:rhs-steady-state}
\end{equation}
where $\mathbf{U}$ is the nodal temperature vector.

For the transient heat conduction, AMATRX and RHS are defined as follows:
\begin{equation}
\mathrm{AMATRX}=\mathbf{K}^{t+\Delta t}+\frac{\mathbf{M}^{t+\Delta t}}{\Delta t}, \label{eq:amatrx-tranisent-state}
\end{equation}
\begin{equation}
    \mathrm{RHS}=-\mathbf{K}^{t+\Delta t}\mathbf{U}^{t+\Delta t}-\frac{\mathbf{M}^{t+\Delta t}}{\Delta t}\Big(\mathbf{U}^{t+\Delta t}-\mathbf{U}^{t}\Big). \label{eq:rhs-tranisent-state}
\end{equation}

\begin{figure}[H]
  \centering
  \includegraphics[width=0.9\textwidth]{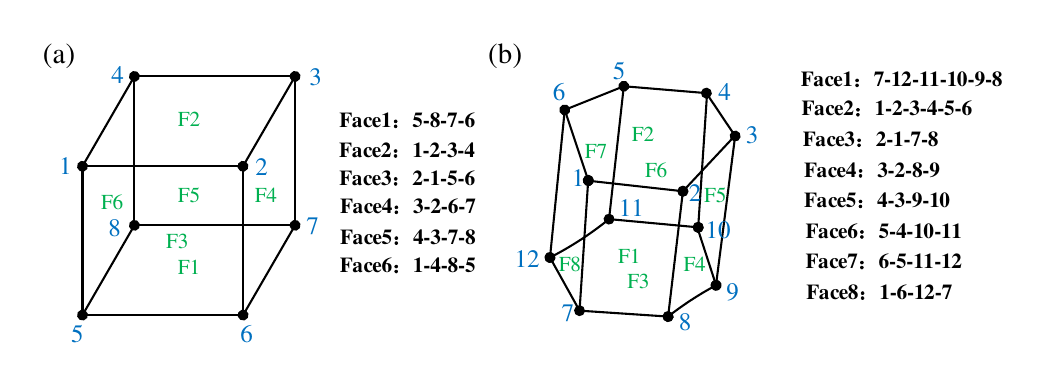}
  \caption{An example of face numbering for polyhedral elements: (a) hexahedral element; (b) polyhedral element.}
  \label{fig:An example of face numbering on polyhedral elements.}
\end{figure}

Moreover, it is necessary to construct the element geometry using nodal coordinates and node connectivity information. The element is first constructed to define the numbering and connectivity of the element nodes. As illustrated in Figs. \ref{fig:An example of face numbering on polyhedral elements.} (a) and (b), a hexahedral element and an octahedral element are considered. The numbering of the element faces is defined, followed by the connectivity information of the nodes on each face.

\begin{algorithm}
	\caption{Solving the 3D heat conduction problems using the PSBFEM}\label{alg:1}
    \textbf{Input:} Node and element information, material properties, and nodal temperature $\mathbf{u}_t$\\
    \textbf{Output:} Nodal temperature $\mathbf{u}_{t+1}$
	\begin{algorithmic}[1]
        \While{ABAQUS not converged}
        \State{Solve nodal temperature $\mathbf{u}^k_{t+1}$}
		\For {1 to AllEle}  \quad \quad  \textit{! Traverse all elements}
         \State{Construct the geometry of elements}
         \State{Solve the the centre of elements}
         \State{Solve the coefficient matrices $\mathbf{E^e_0}$,$\mathbf{E^e_1}$,$\mathbf{E^e_2}$,and $\mathbf{M^e_0}$ using Eq. (\ref{eq:E0}) $\sim$ Eq. (\ref{eq:M0})}
         \State{Construct $\mathbf{Z_p}$ using Eq. (\ref{eq:zp})}
         \State{Solve the $\mathbf{Z_p}$ by eigenvalue decomposition using Eq. (\ref{eq:eigen decomp})}
         \State{solve the stiffness $\mathbf{K}$ and mass matrix $\mathbf{M}$ using Eq. (\ref{eq:K}) and Eq. (\ref{eq:M})}
         \If{lflags(1)=31}
         \State{Update stiffness matrix $\mathbf{AMATRX}$ and residual force vetor $\mathbf{RHS}$ using Eqs. (\ref{eq:amatrx-steady-state}) and (\ref{eq:rhs-steady-state})}
         \EndIf
         \If{lflags(1)=32 or lflags(1)=33}
         \State{Update stiffness matrix $\mathbf{AMATRX}$ and residual force vetor $\mathbf{RHS}$ using Eqs. (\ref{eq:amatrx-tranisent-state}) and (\ref{eq:rhs-tranisent-state})}
         \EndIf
		\EndFor
        \State{\textbf{update} $k=k+1$}
        \State{Solve nodal temperature $\mathbf{u}_{t+1}$=$\mathbf{u}^k_{t+1}$}
        \EndWhile
	\end{algorithmic} 
\end{algorithm} 

\color{blue}
\subsection{Acceleration technology based on the octree parent element}
In large-scale numerical simulations, models are generally made up of numerous elements, each requiring individual computation. Improving the computational efficiency of a single element can therefore lead to faster overall model computation. A key benefit of employing octree discretization is that cubic elements, especially those without hanging nodes, constitute a significant portion of the mesh elements \cite{chen2018efficient,zou2019approach}. 

\begin{figure}[H]
  \centering
  \includegraphics[width=1.0\textwidth]{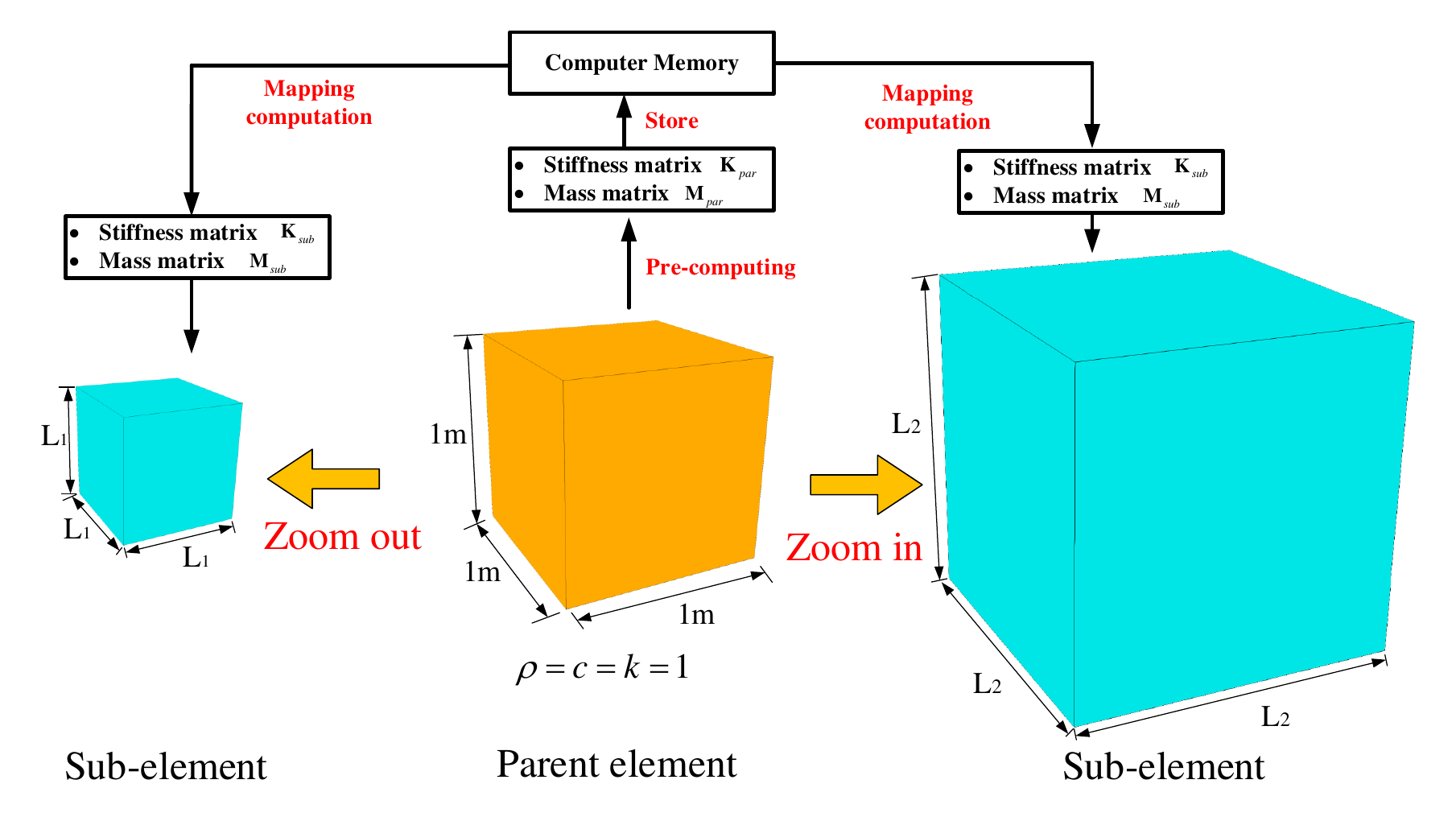}
  \caption{\textcolor{blue}{Mapping relationship between the parent element and sub-element.}}
  \label{fig:parent_and_sub}
\end{figure}

We define a cubic parent element (without hanging nodes) whose geometric and material parameters are set to 1, as shown in Fig. \ref{fig:parent_and_sub}. The mapping relationship between the parent element and the sub-element can be expressed as
\begin{equation}
    \mathbf{K}_{sub}=kL\mathbf{K}_{par},\label{eq:map1}
\end{equation}
\begin{equation}
    \mathbf{M}_{sub}=\rho c L^3\mathbf{M}_{par},\label{eq:map2}
\end{equation}
where $\mathbf{K}_{sub}$ and $\mathbf{M}_{sub}$ are the stiffness matrix and mass matrix of the sub-element, respectively, and $\mathbf{K}_{par}$ and $\mathbf{M}_{par}$ are the stiffness matrix and mass matrix of the parent element, respectively. $k$, $\rho$, $c$, and $L$ represent the thermal conductivity, density, heat capacity, and length of the sub-element, respectively.

Before the analysis, the stiffness matrix $\mathbf{K}_{par}$ and mass matrix $\mathbf{M}_{par}$ of the parent element are computed and stored in the computer memory. During the analysis, when the element is a cubic element, the stiffness and mass matrices of the sub-elements are computed based on the mapping relationships in Eqs. (\ref{eq:map1}) and (\ref{eq:map2}). Only the material and geometric coefficients need to be multiplied, eliminating the need to apply the SBFEM method to compute the stiffness and mass matrices. As a result, this approach significantly reduces the computational cost.

\color{black}
\subsection{Defining the element of the UEL}
\begin{figure}[H]
  \centering
  \includegraphics[width=0.8\textwidth]{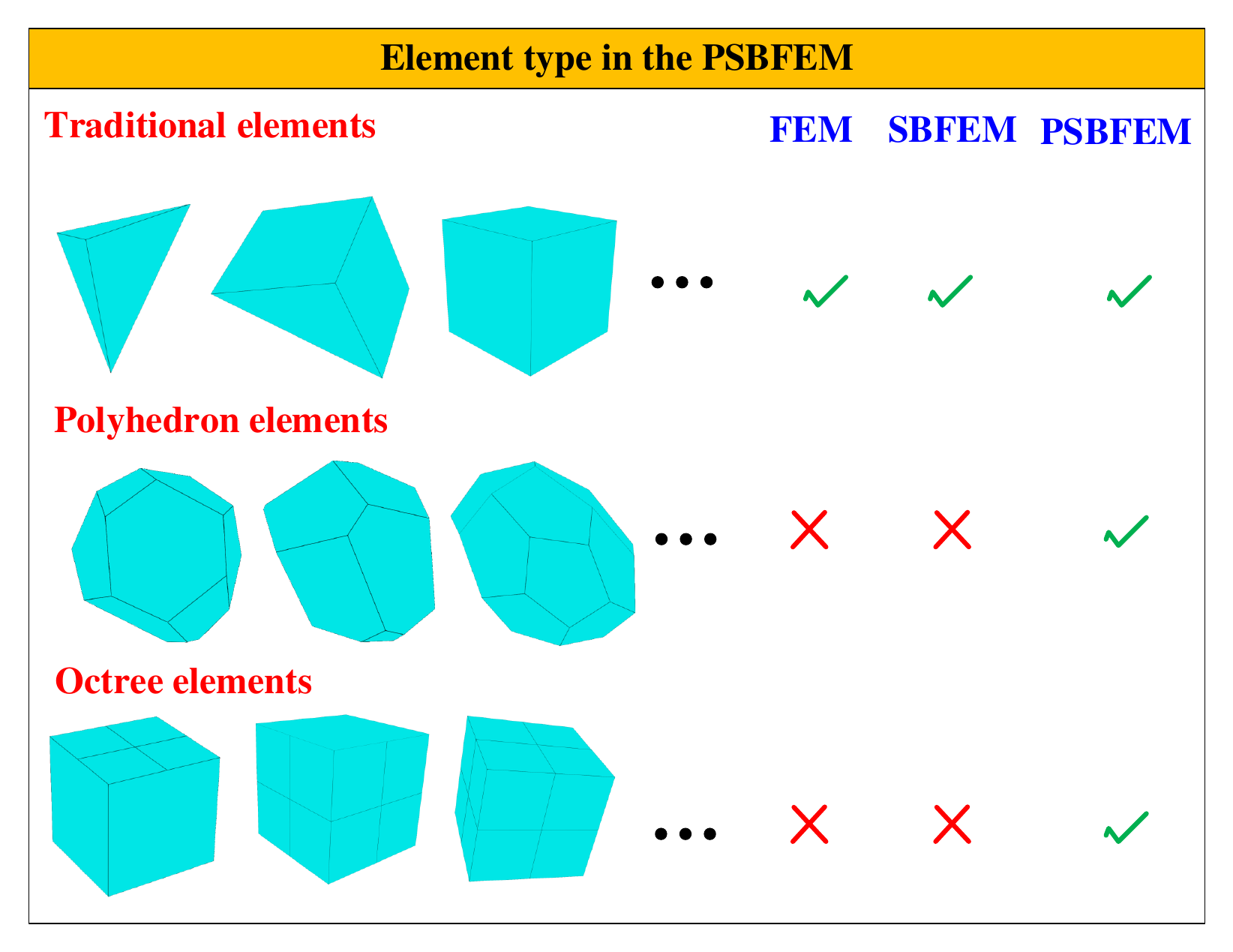}
  \caption{The supported elements of the PSBFEM.}
  \label{fig:supported element}
\end{figure}

In 3D PSBFEM analysis, the geometry of the elements is constructed based on the nodes. Fig. \ref{fig:An example of face numbering on polyhedral elements.} illustrates an example of element construction. Prior to the analysis, it is essential to establish the geometric relationships of the element surfaces according to their numbering, ensuring that the Jacobian matrix remains negative. This requires that the direction of each surface points outward from the element. Fig. \ref{fig:supported element} shows the supported element types in the PSBFEM. In the PSBFEM, traditional elements such as tetrahedral, wedge, and hexahedral elements can be used. Additionally, the PSBFEM also supports more complex element types, including polyhedral and octree elements. Thus, the PSBFEM serves as an effective tool for solving complex geometries in numerical analysis

The ABAQUS input file typically provides a comprehensive representation of the numerical model, detailing aspects such as nodes, elements, degrees of freedom, and materials. As shown in Listing 1, the hexahedral element (U8) is defined as follows: Lines 1 to 4 are used to define the pentagonal element (U8). Line 1 specifies the element type, number of nodes, number of element properties, and number of degrees of freedom for each node; Line 2 defines the active degrees of freedom for temperature; Lines 3 and 4 assign the element set labeled \textbf{Hexahedron}. Similarly, other elements can be defined using the same procedure.

\begin{lstlisting}[caption={Input file of the polyhedral element in UEL.}]
*USER ELEMENT,NODES=8,TYPE=U8,PROPERTIES=5,COORDINATES=3,VARIABLES=8
11
*Element,TYPE=U8,ELSET=Hexahedron
1, 5, 6, 8, 7, 1, 2, 4, 3
*USER ELEMENT,NODES=12,TYPE=U12,PROPERTIES=5,COORDINATES=3,VARIABLES=12
11
*Element,TYPE=U12,ELSET=Octahedron
2, 7, 8, 9, 10, 11, 12, 1, 2, 3, 4, 5, 6
\end{lstlisting}

\section{Numerical examples}
\label{sec:7}
In this section, several benchmark problems were presented to demonstrate the convergence and accuracy of the proposed framework for heat conduction analysis. Moreover, the results obtained from the PSBFEM were compared with those from the FEM, with the FEM analysis performed using the commercial software ABAQUS. The comparisons were conducted on a system equipped with an Intel Core i5-10210U CPU (1.60 GHz) and 8.0 GB of RAM. To validate the proposed approach, the relative errors $\mathbf{e}_{L_2}$ in temperature were examined as follows:
\begin{equation}
    \mathbf{e}_{L_2}=\frac{\left\|\mathbf{T}_{num}-\mathbf{T}_{ref}\right\|_{L_2}}{\parallel\mathbf{T}_{ref}\parallel_{L_2}},
\end{equation}
where $\mathbf{T}_{num}$ is the numerical results and $\mathbf{T}_{ref}$ is the reference solutions.

\subsection{Patch test}
\begin{figure}[H]
  \centering
  \includegraphics[width=0.8\textwidth]{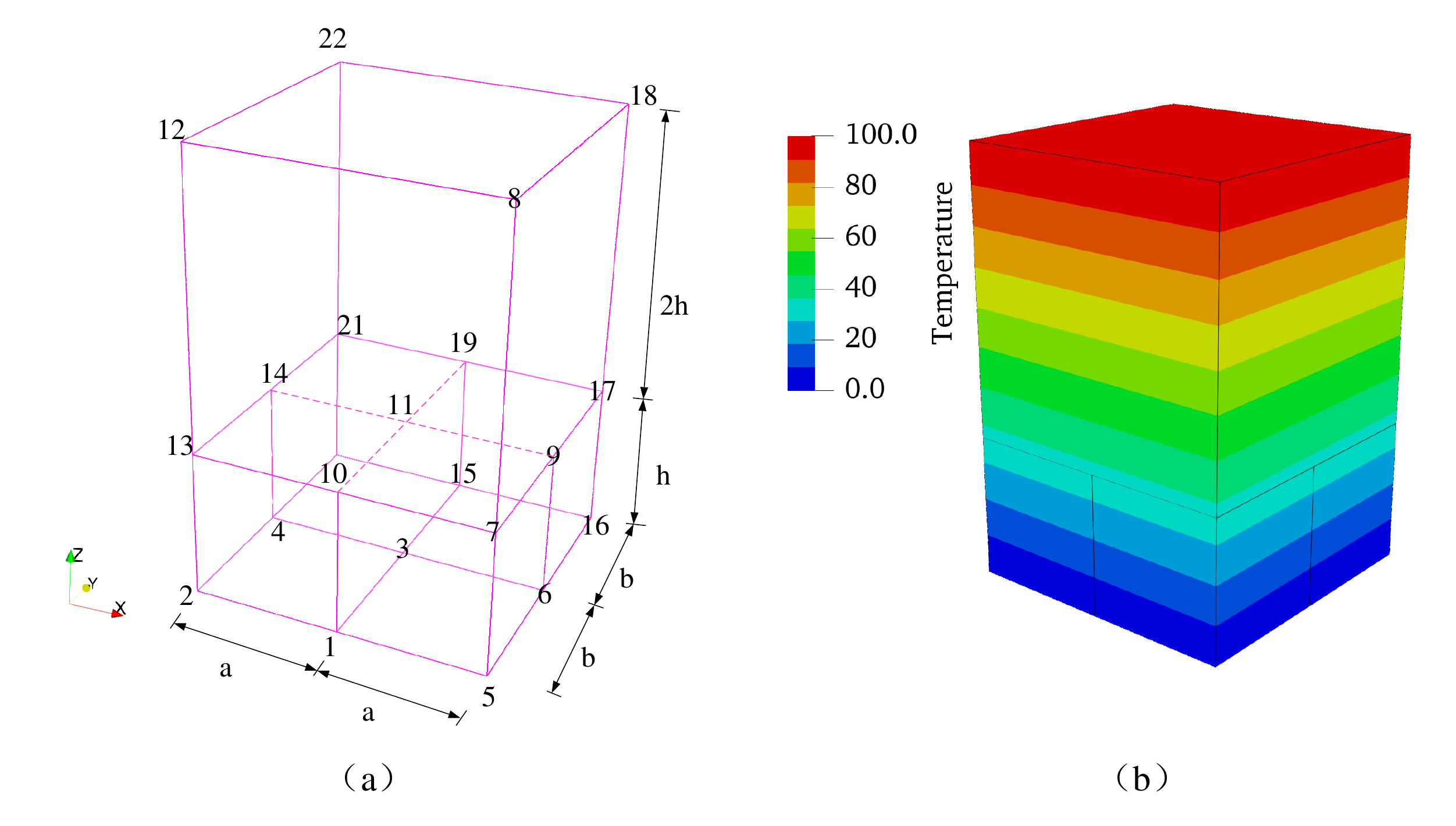}
  \caption{Geometry model and temperature contour of the patch test (Unit: $^{\circ}\mathrm{C}$).}
  \label{fig:patch_geo_contour}
\end{figure}

To verify that the proposed method satisfies the convergence requirement, it must pass the standard patch test \cite{yaOpen2021,zhangFast2019}. As shown in Fig. \ref{fig:patch_geo_contour} (a), the patch test was conducted using a quadrangular prism (\(a=1\)  m, \(b=1\)  m, \(h=1\) m). The material constants were a thermal conductivity of 1.0 W/m/\(^\circ\mathrm{C}\) and a volumetric heat capacity of \(\rho c = 1.0 \mathrm{J}/(\mathrm{m}^3 \cdot ^\circ\mathrm{C})\). Temperatures of \(T=100^\circ\mathrm{C}\) and \(T=0^\circ\mathrm{C}\) were prescribed on the top surface (\(z=3\) m) and bottom surface (\(z=0\) m), respectively. The temperature distribution was shown in Fig. \ref{fig:patch_geo_contour} (b). Tab. \ref{tab:patch} presented the relative errors for nodal temperatures calculated at all nodes and compared to the analytical solution. The maximum relative error of the proposed method in temperature compared to the analytical solution was \(1.53 \times 10^{-4}\). Thus, the proposed method passed the patch test with sufficient accuracy.

\begin{table}[H]
\centering
\caption{Maximum relative error of the nodal temperature.}
\begin{tabular}{@{}lll@{}}
\toprule
 Proposed method & Analytical solution & Relative error\\ \midrule
 33.3282 & 33.3333 & $1.53 \times 10^{-4}$\\ \bottomrule
\end{tabular}
\label{tab:patch}
\end{table}

\subsection{3D heat conduction beam}
\begin{figure}[H]
  \centering
  \includegraphics[width=0.8\textwidth]{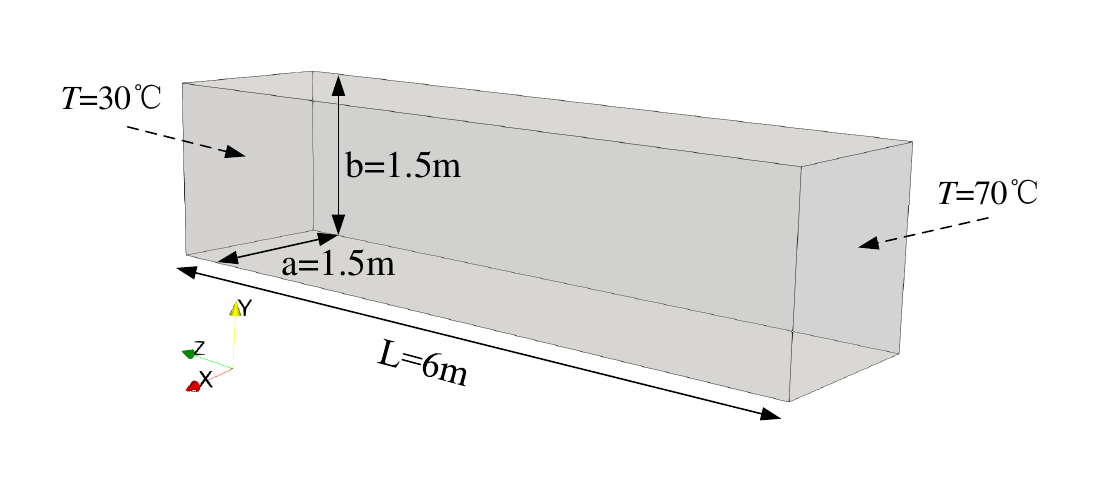}
  \caption{Geometry and boundary conditions for a 3D heat conduction beam.}
  \label{fig:ex01_geo}
\end{figure}

In this example, we examined a 3D heat conduction problem, as illustrated in Fig. \ref{fig:ex01_geo}. The beam had a length of 6 m, with a width and height of 1.5 m. A temperature of \(T=70^{\circ}\mathrm{C}\) was applied at the right boundary (\(z=0\) m), while a temperature of \(T=30^{\circ}\mathrm{C}\) was imposed at the left boundary (\(z=6\) m). The material had a thermal conductivity of 1.0 W/m/\(^\circ\mathrm{C}\) and a volumetric heat capacity of \(\rho c = 1.0 \, \mathrm{J}/(\mathrm{m}^3 \cdot ^\circ\mathrm{C})\).

\begin{figure}[H]
  \centering
  \includegraphics[width=0.8\textwidth]{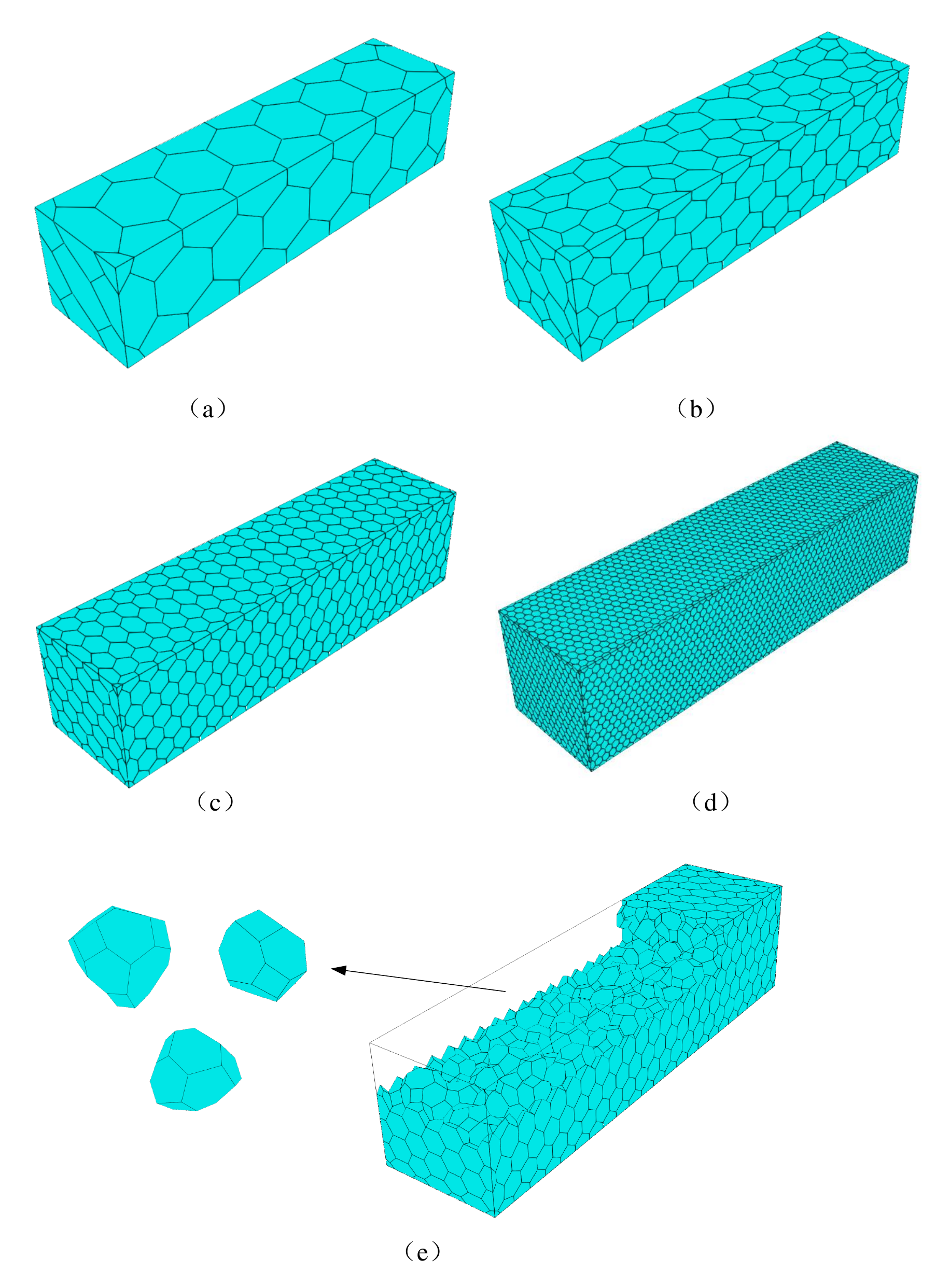}
  \caption{Four levels of polyhedral mesh refinement: (a) 1.0 m mesh size; (b) 0.5 m mesh size; (c) 0.25 m mesh size; (d) 0.1 m mesh size; (e) internal element of the 0.25 m mesh size.}
  \label{fig:ex01_polymesh}
\end{figure}

\begin{figure}[H]
  \centering
  \includegraphics[width=1.0\textwidth]{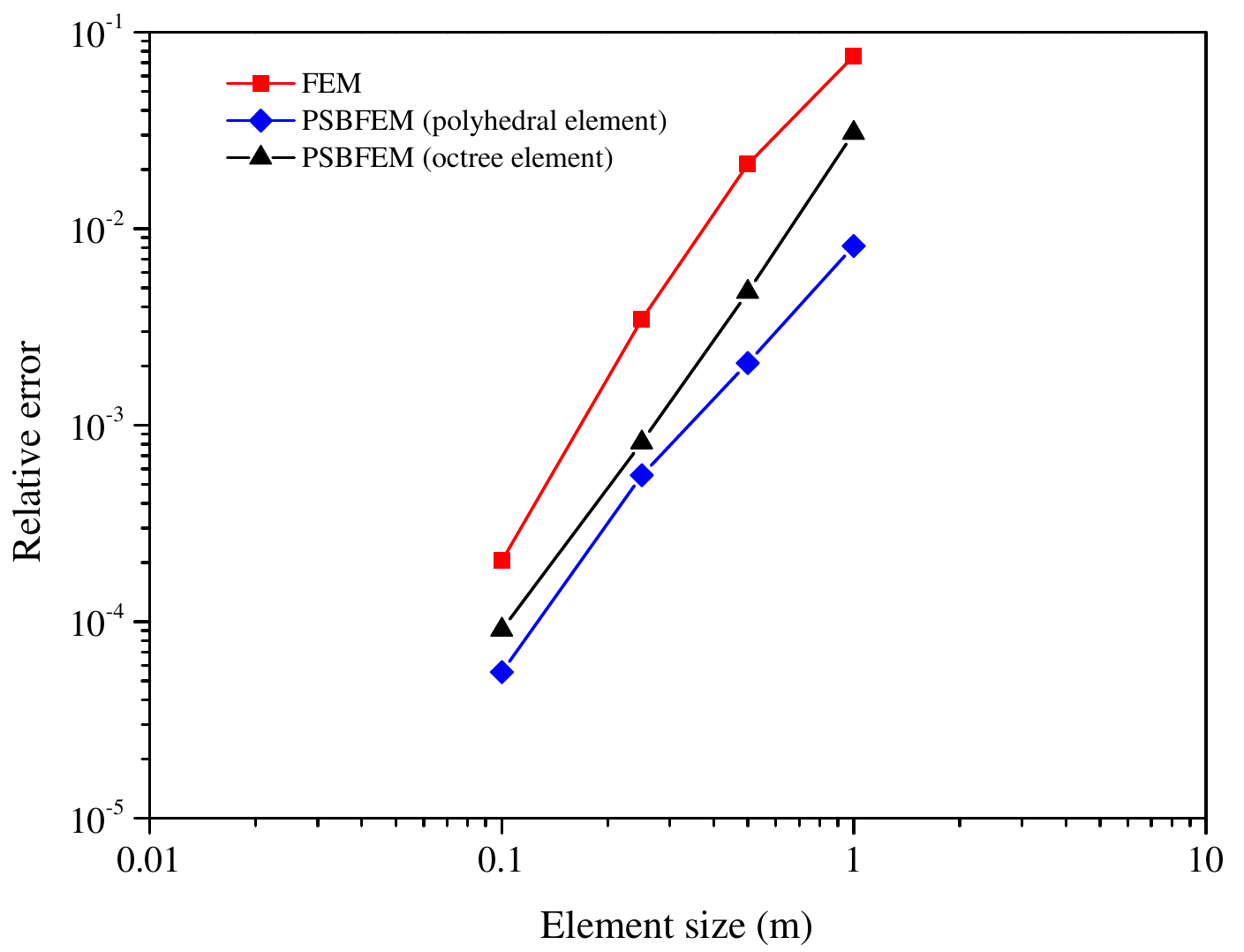}
  \caption{\color{blue}{Convergence of the relative errors in the temperature.}}
  \label{fig:ex01_conv}
\end{figure}

\begin{figure}[H]
  \centering
  \includegraphics[width=0.8\textwidth]{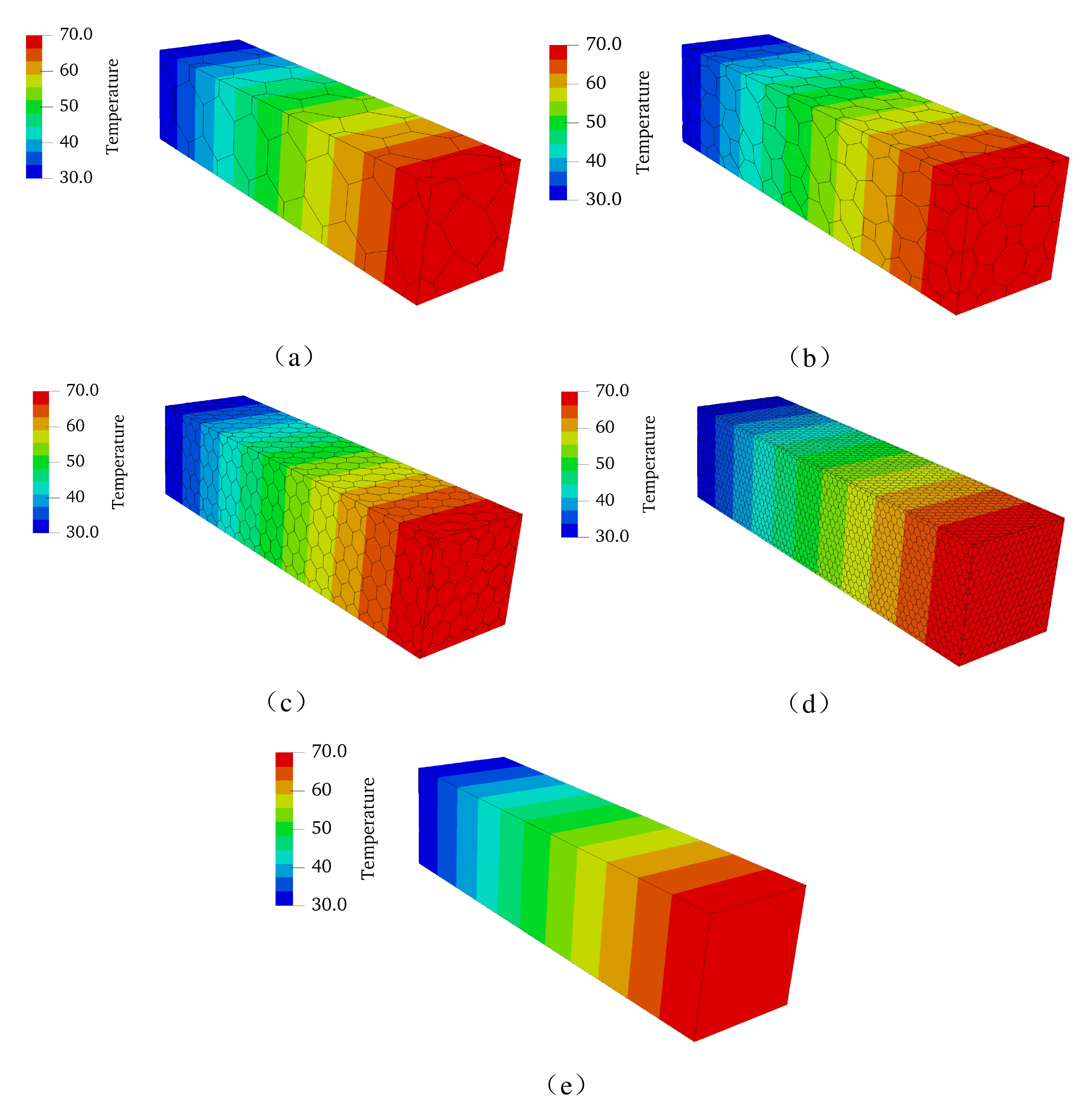}
  \caption{Temperature distribution of the 3D heat conduction beam using polyhedral mesh refinement: (a) 1.0 m mesh size; (b) 0.5 m mesh size; (c) 0.25 m mesh size; (d) 0.1 m mesh size; (e) analytical solution (Unit: $^{\circ}\mathrm{C}$).}
  \label{fig:ex01_polyresult}
\end{figure}

\textcolor{blue}{The domain was discretized using hexahedral, polyhedral and octree elements}. A convergence study was conducted with all three types of mesh. We refined the meshes successively following the sequence \(h = 1.0\) m, \(0.5\) m, \(0.25\) m, and \(0.1 \, \mathrm{m}\), as shown in Fig. \ref{fig:ex01_polymesh}. Moreover, Fig. \ref{fig:ex01_polymesh} (e) presented a cross-sectional view of a mesh layer along with several representative polyhedral elements. Fig. \ref{fig:ex01_conv} showed that the PSBFEM and FEM exhibit a good convergence rate, with the accuracy of both methods increasing as the mesh was refined. \textcolor{blue}{When the mesh size was \(0.1 \, \mathrm{m}\), the relative errors for the FEM, PSBFEM (polyhedral element) and PSBFEM (octree element) were \(2.1 \times 10^{-4}\), \(5.5 \times 10^{-5}\) and \(9.1 \times 10^{-5}\) respectively.} The relative errors of the PSBFEM were smaller than those of the FEM for the same element size. \textcolor{blue}{The relative errors for the polyhedral elements in PSBFEM were smaller than those for the octree elements at the same element size. This was mainly due to the fact that polyhedral elements have more nodes at the same element size.} Moreover, Fig. \ref{fig:ex01_polyresult} demonstrated a significant concurrence in the temperature distribution between the PSBFEM and the analytical solution with refined meshes. 

\begin{figure}[H]
  \centering
  \includegraphics[width=0.8\textwidth]{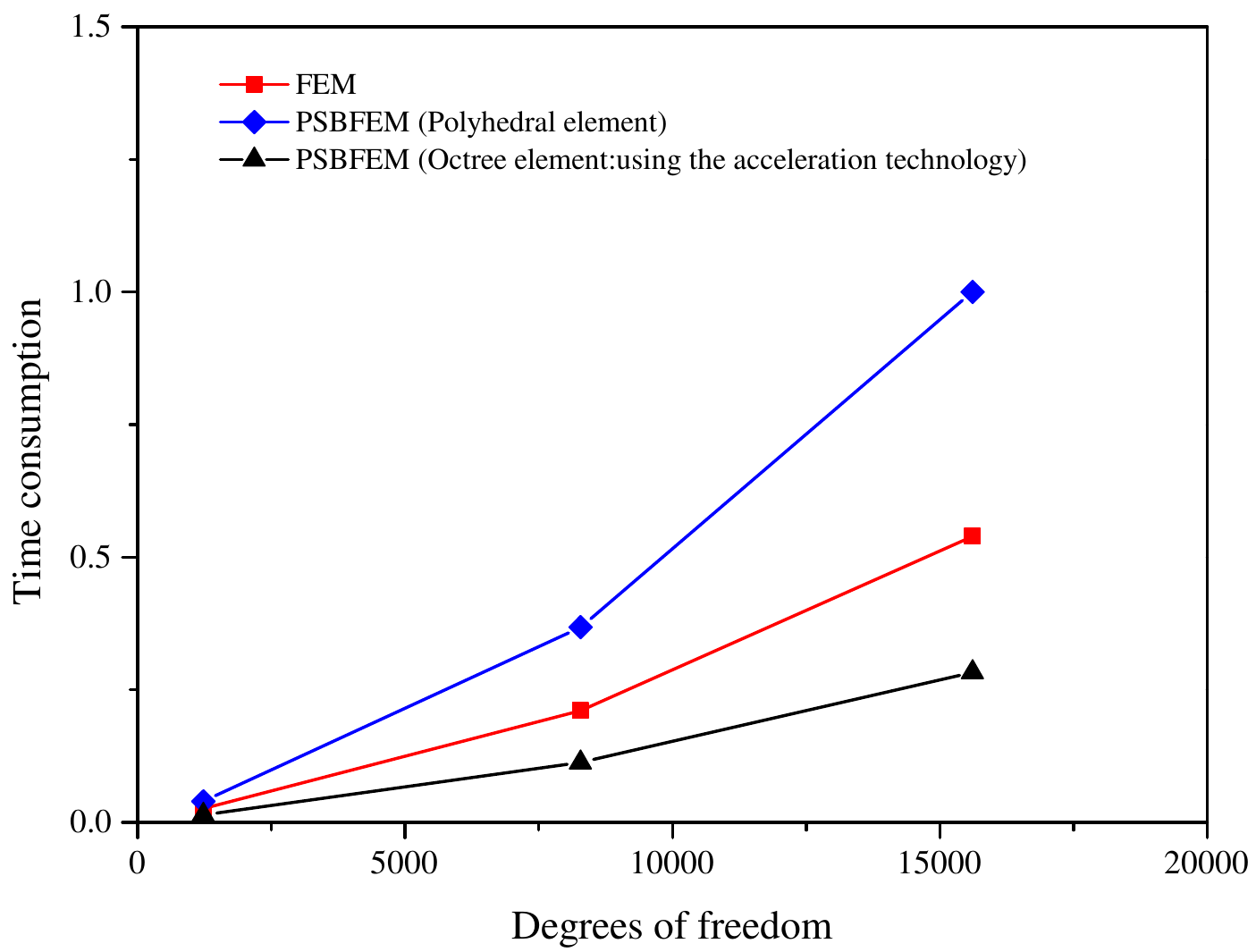}
  \caption{\textcolor{blue}{Time consumption comparison for the steady-state heat conduction analysis.}}
  \label{fig:ex01_time}
\end{figure}

\textcolor{blue}{Additionally, Fig. \ref{fig:ex01_time} compared the computational cost of PSBFEM and FEM under the same degrees of freedom. To ensure a fair comparison, an automatic incrementation scheme was adopted to minimize the impact of increment size on solving time. The computation time was normalized for comparison. The computational cost of PSBFEM using the polyhedral element was slightly higher than that of the FEM, averaging approximately 1.70 times that of the FEM. This increases in computational cost can be attributed to the semi-analytical nature of SBFEM, which inherently involved more computations than FEM \cite{yePSBFEM2021,yu2021scaled}. In addition, by using acceleration technology, the computational time of PSBFEM with the octree mesh was smaller than that of FEM, averaging 0.53 times that of the FEM.}

\begin{figure}[H]
  \centering
  \includegraphics[width=0.9\textwidth]{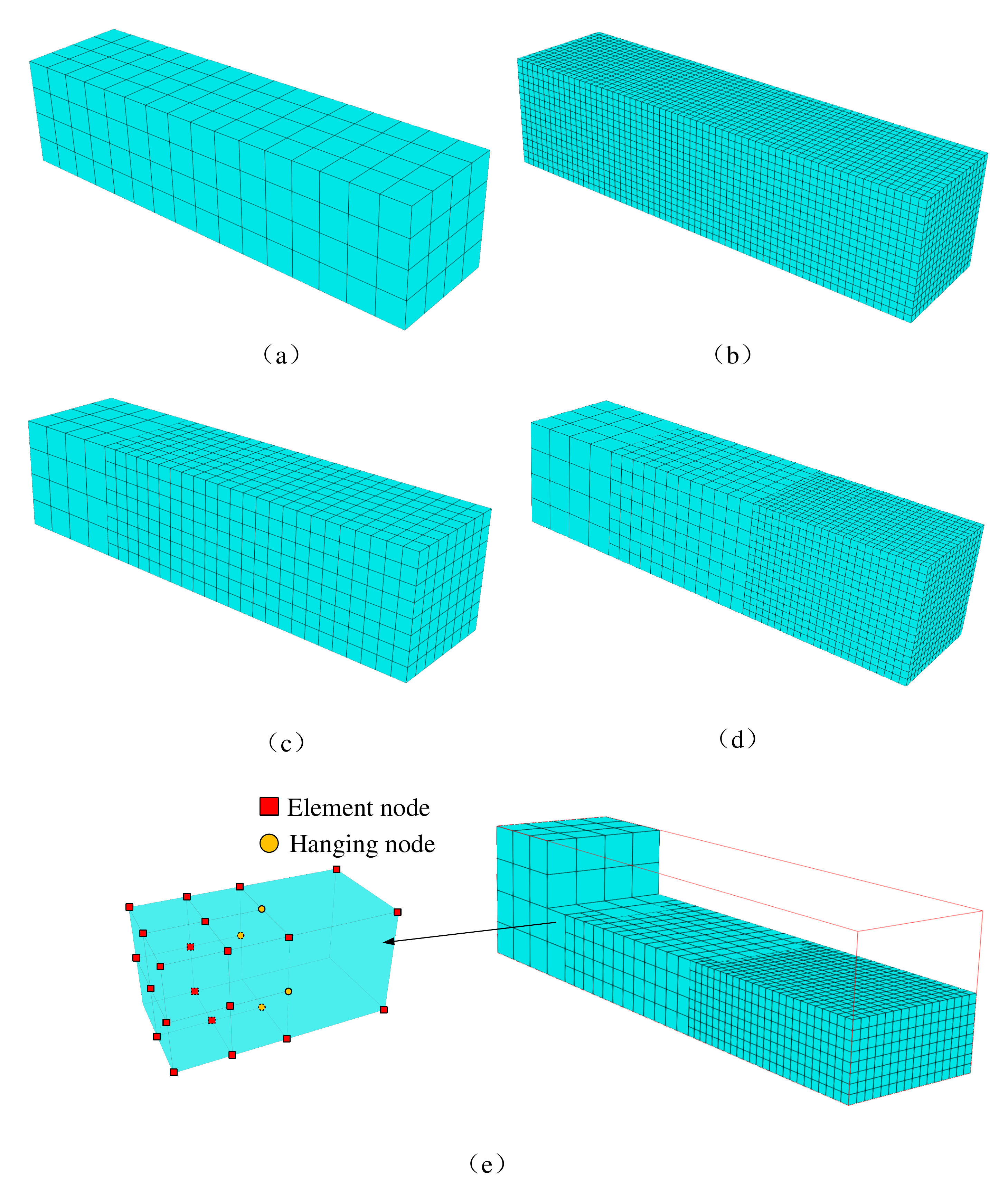}
  \caption{\textcolor{blue}{Local mesh refinement model of the 3D heat conduction beam using the octree mesh: (a) coarse mesh; (b) fine mesh; (c) local mesh refinement 1; (d) local mesh refinement 2; (e) internal element of local mesh refinement 2.}}
  \label{fig:ex01_octreeMesh}
\end{figure}

\begin{figure}[H]
  \centering
  \includegraphics[width=0.9\textwidth]{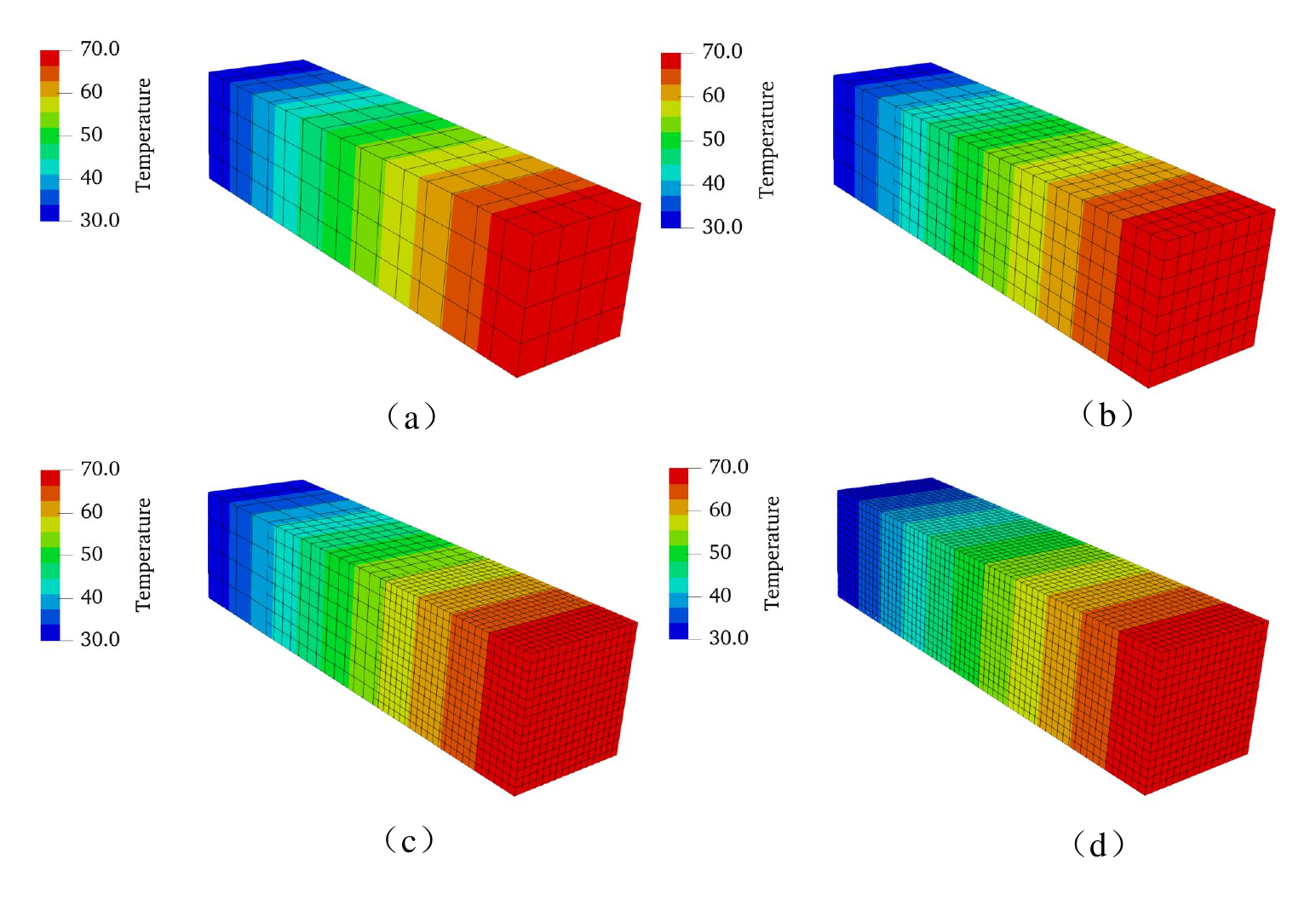}
  \caption{Temperature distribution of the 3D conduction beam using the octree mesh: (a) temperature contour with coarse mesh; (b) temperature contour with local mesh refinement 1; (c) temperature contour with local mesh refinement 2; (d) temperature contour with fine mesh; (Unit: $^{\circ}\mathrm{C}$).}
  \label{fig:ex01_octreeresult}
\end{figure}

\textcolor{blue}{To demonstrate the capability of octree mesh discretization in local mesh refinement, we discretized a 3D heat conduction beam using an octree mesh,}
as shown in Fig. \ref{fig:ex01_octreeMesh}. Figs. \ref{fig:ex01_octreeMesh} (a) and \ref{fig:ex01_octreeMesh} (b) represent the coarse and fine meshes, respectively, while Figs. \ref{fig:ex01_octreeMesh} (c) and \ref{fig:ex01_octreeMesh} (d) depict the local mesh refinement. In addition, Fig. \ref{fig:ex01_octreeMesh} (e) presented a cross-section of local mesh refinement 2. Tab. \ref{tab:ex01_t1} illustrated the mesh characteristics and relative errors for the local mesh refinement. The relative errors for local mesh refinement 2 and the fine mesh are both on the order of \(10^{-4}\), demonstrating a high level of computational accuracy. However, the computational time for local mesh refinement 2 was lower than that of the fine mesh. Thus, through local mesh refinement, refining the mesh in specific regions can yield more accurate results while reducing computational time. \textcolor{blue}{As shown in Tab. \ref{tab:ex01_t1}, with the acceleration technology, the computational time of PBSFEM was significantly reduced, with an average reduction of 3 times in the computational time.} Moreover, Fig. \ref{fig:ex01_octreeresult} also showed a close match in the temperature distribution between the analytical solution and the PSBFEM solution using the octree mesh for discretization.
 
\begin{table}[H]
\centering
\normalsize
\caption{Mesh characteristics, relative errors and time consumption.}
\resizebox{\textwidth}{!}{
\begin{tabular}{ccccccc}
\toprule
Mesh type  & Elements & Nodes & Faces &  Relative error & CPU time (s)&\textcolor{blue}{CPU time (s)*}\\ 
\midrule
Coarse mesh            & 256       & 425       & 1536       & 1.78$\times10^{-2}$  & 0.2  &\textcolor{blue}{0.1}\\ 
Fine mesh              & 8432      & 12141     & 53328     & 6.33$\times10^{-4}$   & 9.2   &\textcolor{blue}{3.1}\\ 
\textcolor{blue}{Local mesh refinement 1}     & 1572      & 2108      & 9504       & 2.80$\times10^{-3}$  & 1.5   &\textcolor{blue}{0.4}\\ 
\textcolor{blue}{Local mesh refinement 2}     & 6696      & 8063      & 40512     & 8.49$\times10^{-4}$   & 6.2   &\textcolor{blue}{1.8}\\ 
\bottomrule 
\label{tab:ex01_t1}
\end{tabular}}
\vspace{-2em}
\begin{flushleft}
\textcolor{blue}{\small Note: * denotes the use of acceleration technology based on the parent element.}
\end{flushleft}
\end{table}

\subsection{Steady-state heat conduction}
\begin{figure}[H]
  \centering
  \includegraphics[width=0.9\textwidth]{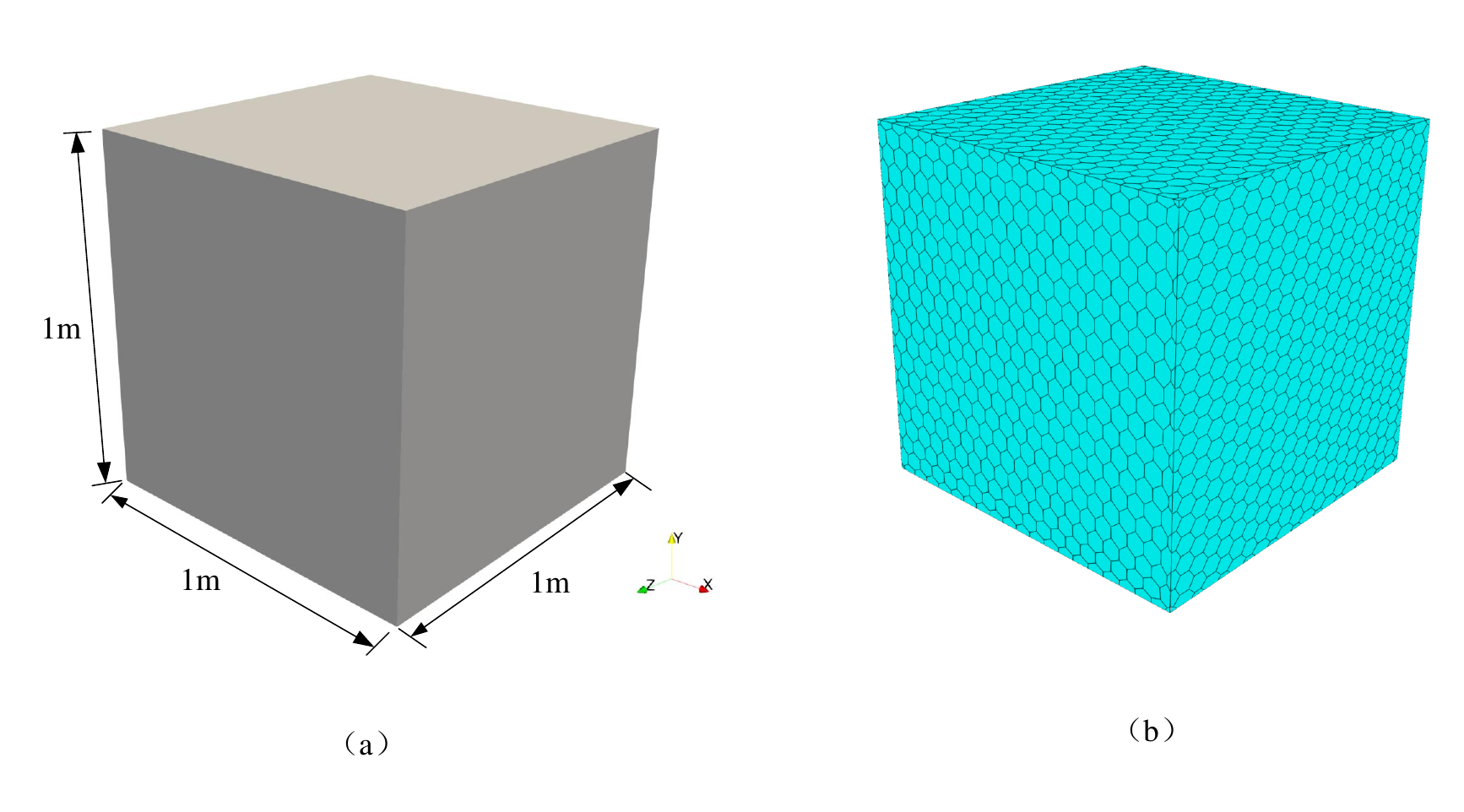}
  \caption{A cube for the steady-state heat conduction problem: (a) the geometry model; (b) the mesh with a size of 0.05 m.}
  \label{fig:ex02geo_mesh}
\end{figure}

In this example, we considered a classical steady-state heat conduction problem, as shown in Fig. \ref{fig:ex02geo_mesh}. The length, width, and height of the cube were all 1 m. The material had a thermal conductivity of 1.0 W/m/\(^\circ\mathrm{C}\) and a volumetric heat capacity of \(\rho c = 1.0 \, \mathrm{J}/(\mathrm{m}^3 \cdot ^\circ\mathrm{C})\). The boundary conditions of the cube were shown in Eq. (\ref{eq:ex02_bc}).

\begin{equation}
    \begin{aligned}
&(x,y,z)\in plane\perp(x=0),\quad T=0 \\
&(x,y,z)\in plane\perp(y=0),\quad T=0 \\
&(x,y,z)\in plane\perp(z=0),\quad T=0 \\
&(x,y,z)\in plane\perp(x=1),\quad T=0 \\
&(x,y,z)\in plane\perp(y=1),\quad T=\sin(\pi x)\sin(\pi z) \\
&(x,y,z)\in plane\perp(z=1),\quad T=0
\end{aligned}.
\label{eq:ex02_bc}
\end{equation}

A convergence study was performed using mesh \textit{h}-refinement, with successive refinements applied to mesh element sizes of 0.25 m, 0.1 m, 0.05 m, and 0.025 m. Fig. \ref{fig:ex02geo_mesh} (b) presented the polyhedral mesh with an element size of 0.05 m. The analytical solution for the temperature at point P ($x,y,z$) within the domain was calculated based on the approach outlined in \cite{afrasiabiCont2020}:

\begin{equation}
    T(x,y,z)=\frac{\sinh(\pi y)}{\sinh(\pi)}\sin(\pi x)\sin(\pi z).
\end{equation}

\begin{figure}[H]
  \centering
  \includegraphics[width=0.9\textwidth]{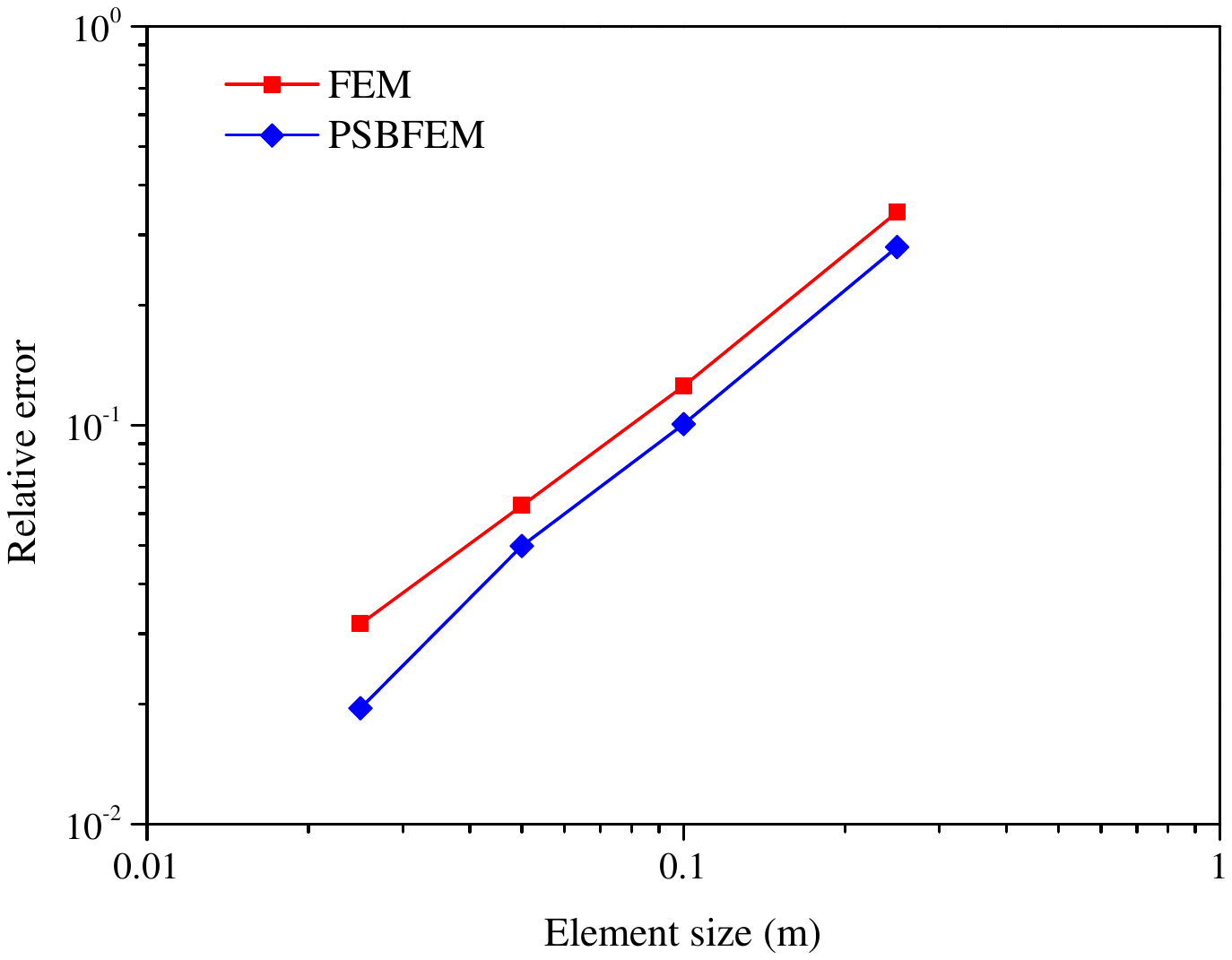}
  \caption{Convergence of the relative error in temperature.}
  \label{fig:ex02_conv}
\end{figure}

\begin{figure}[H]
  \centering
  \includegraphics[width=0.9\textwidth]{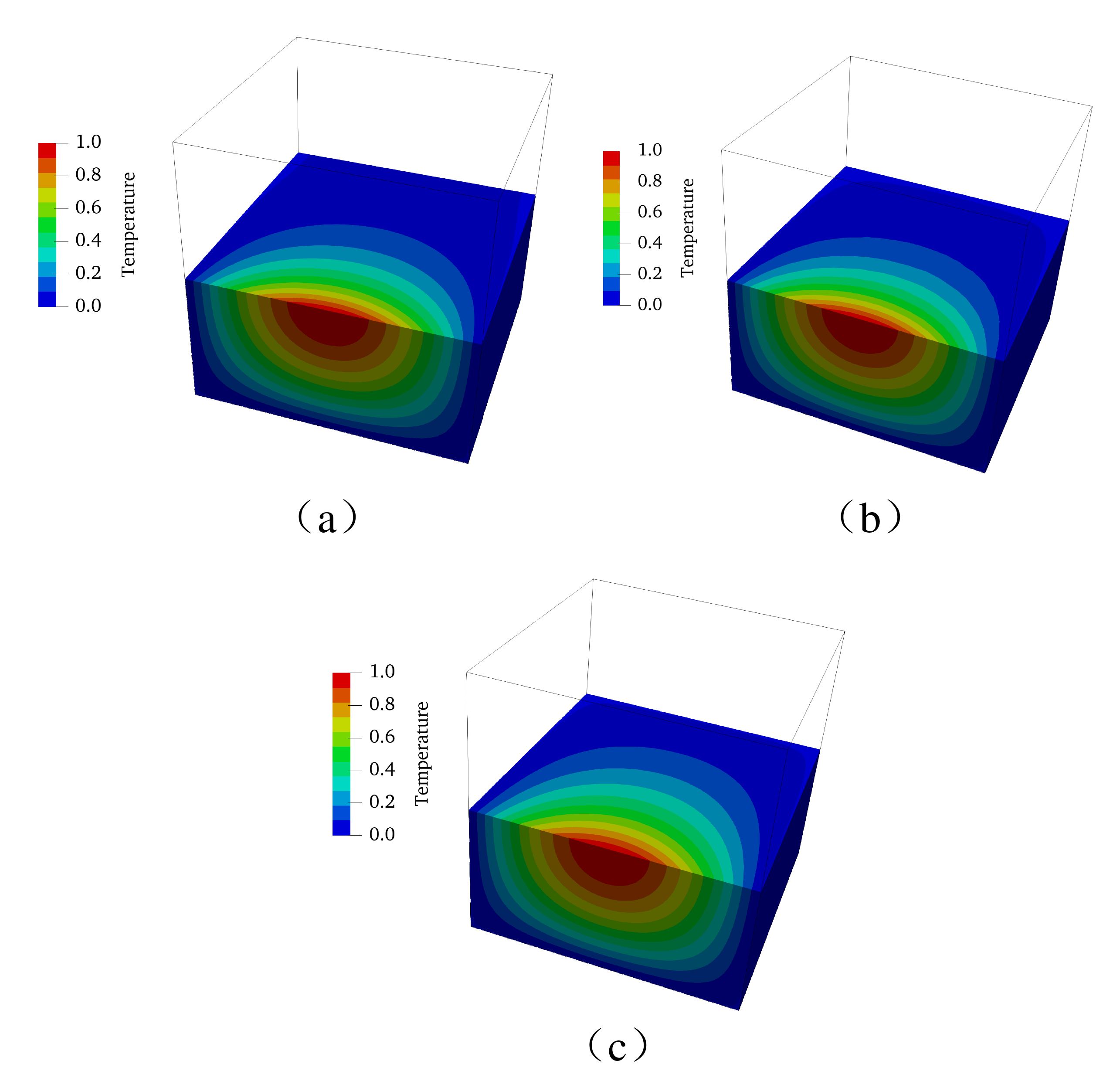}
  \caption{Temperature distribution of a cube for the steady-state heat conduction problem: (a) FEM; (b) PSBFEM; (c) analytical solution; (Unit: $^{\circ}\mathrm{C}$).}
  \label{fig:ex02_contour}
\end{figure}

The convergence of the relative error in temperature with mesh refinement was presented in Fig. \ref{fig:ex02_conv}. It is observed that the PSBFEM converges to the analytical solution at an optimal rate. Furthermore, the PSBFEM produces slightly more accurate results than the FEM for the same element size. Fig. \ref{fig:ex02_contour} illustrated the temperature distribution for both FEM and PSBFEM, showing that the results of both methods closely match the analytical solution.

\subsection{Transient heat conduction}
To validate the proposed PSBFEM for the transient heat conduction problem in three dimensions, we examined a cubic domain \((\pi \times \pi \times \pi)\) with an analytical solution. The initial temperature distribution was defined as \(T_0 = 10 \sin(x) \sin(y) \sin(z)\). Homogeneous Dirichlet boundary conditions were applied on all surfaces of the cube, assuming a constant temperature of zero. The thermal properties of the material were given by a thermal conductivity \(k = 1.0 \, \mathrm{W/m/^\circ\mathrm{C}}\) and a volumetric heat capacity \(\rho c = 1.0 \, \mathrm{J}/(\mathrm{m}^3 \cdot ^\circ\mathrm{C})\). The corresponding analytical solution for the temperature field was expressed as \cite{lin2017transient}:
\begin{equation}
T(x, y, z, t) = 10e^{-3t} \sin(x) \sin(y) \sin(z).
\end{equation}

In both the PSBFEM and FEM simulations, a time step of \(\Delta t = 0.01 \, \mathrm{s}\) was employed, with the total simulation time set to \(t = 1 \, \mathrm{s}\). \textcolor{blue}{ The computational domain was discretized using hexahedral, polyhedral and octree elements.} To assess the accuracy and convergence of the method, a convergence study was conducted through progressive mesh refinement. \textcolor{blue}{The meshes were refined successively in the sequence of $\frac{\pi}{5}$, $\frac{\pi}{10}$, $\frac{\pi}{20}$, and $\frac{\pi}{40}$.}

\begin{figure}[H]
  \centering
  \includegraphics[width=1.0\textwidth]{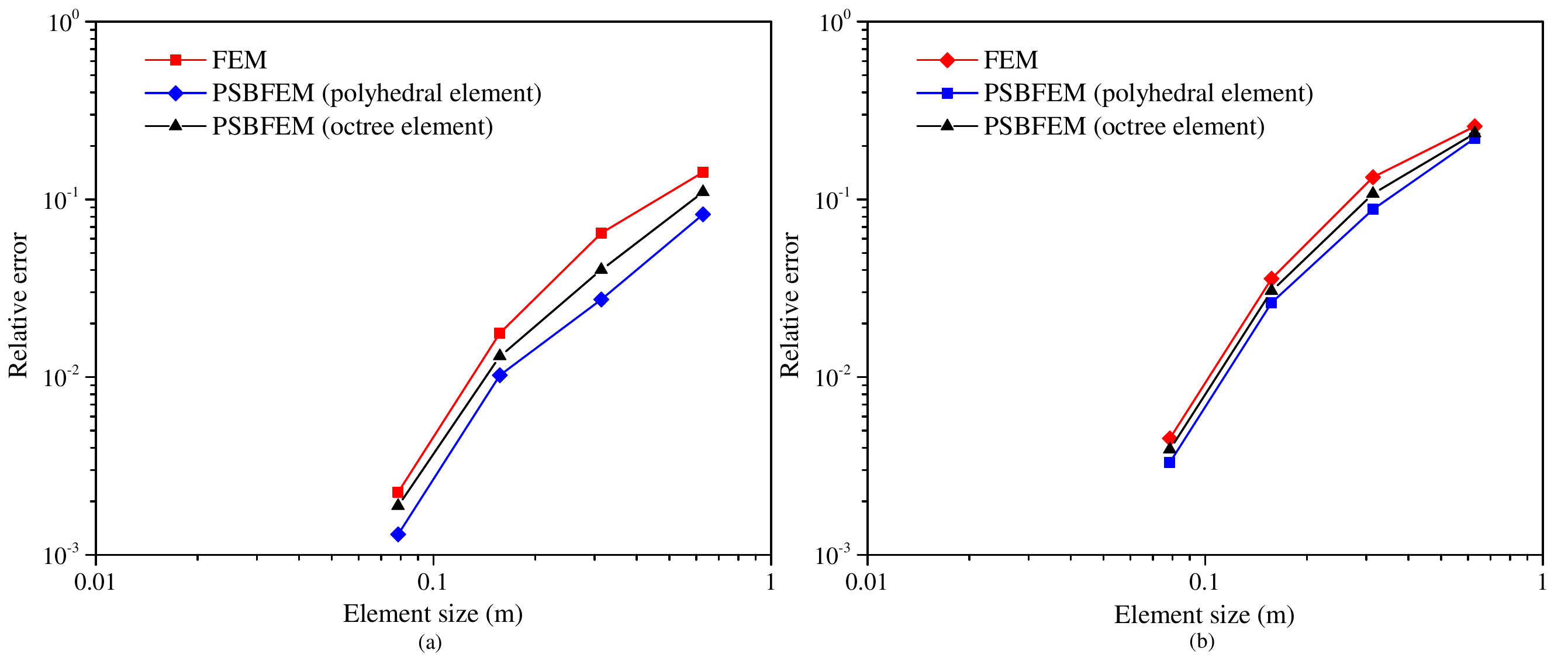}
  \caption{\textcolor{blue}{Convergence of the relative error in temperature at different times: (a) \(t = 0.5 \, \mathrm{s}\); (b) \(t = 1.0 \, \mathrm{s}\).}}
  \label{fig:ex03_conv}
\end{figure}

\begin{figure}[H]
  \centering
  \includegraphics[width=0.7\textwidth]{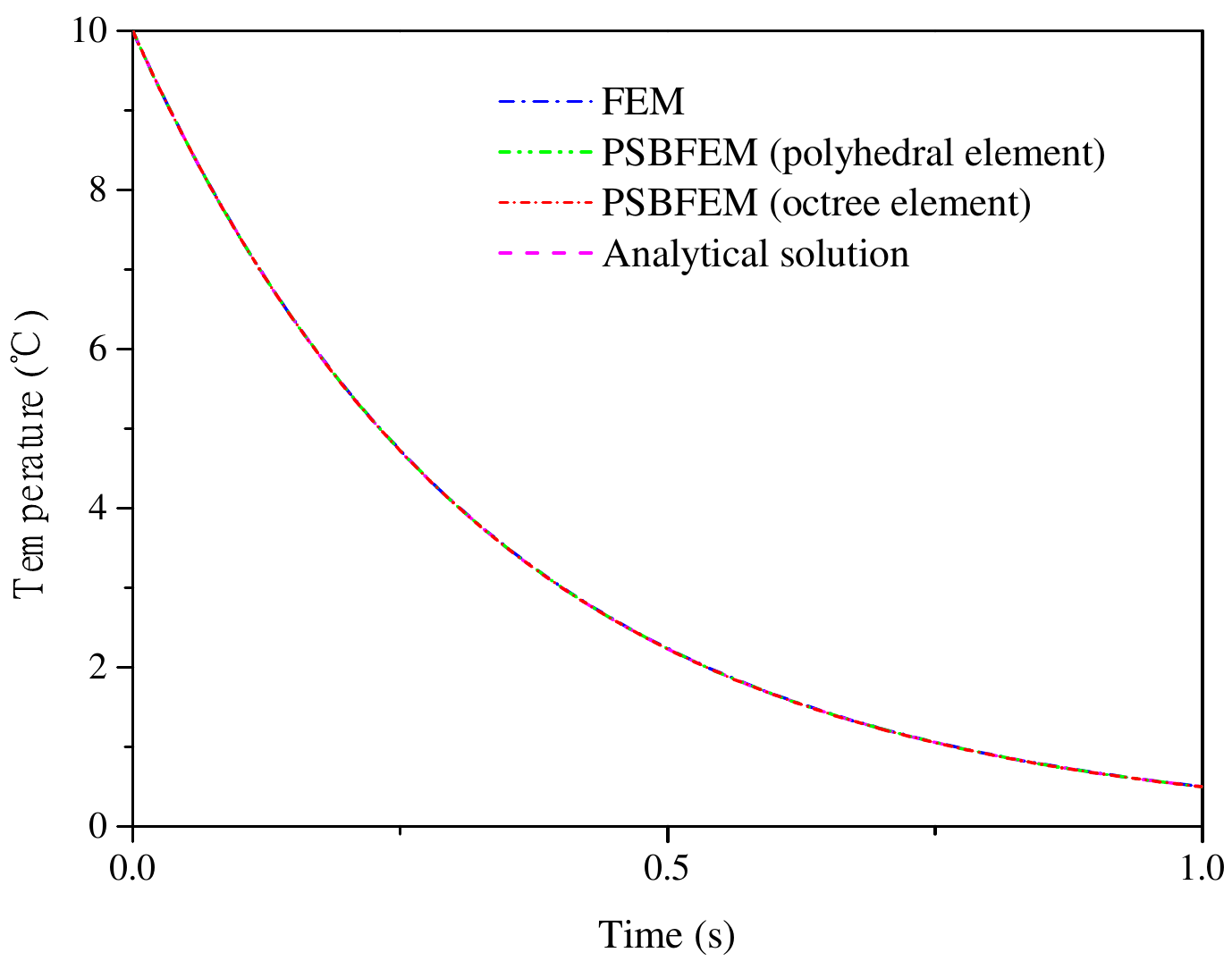}
  \caption{\textcolor{blue}{Comparison of the temperature time history at the monitoring point.}}
  \label{fig:ex03_his}
\end{figure}

\begin{figure}[H]
  \centering
  \includegraphics[width=0.9\textwidth]{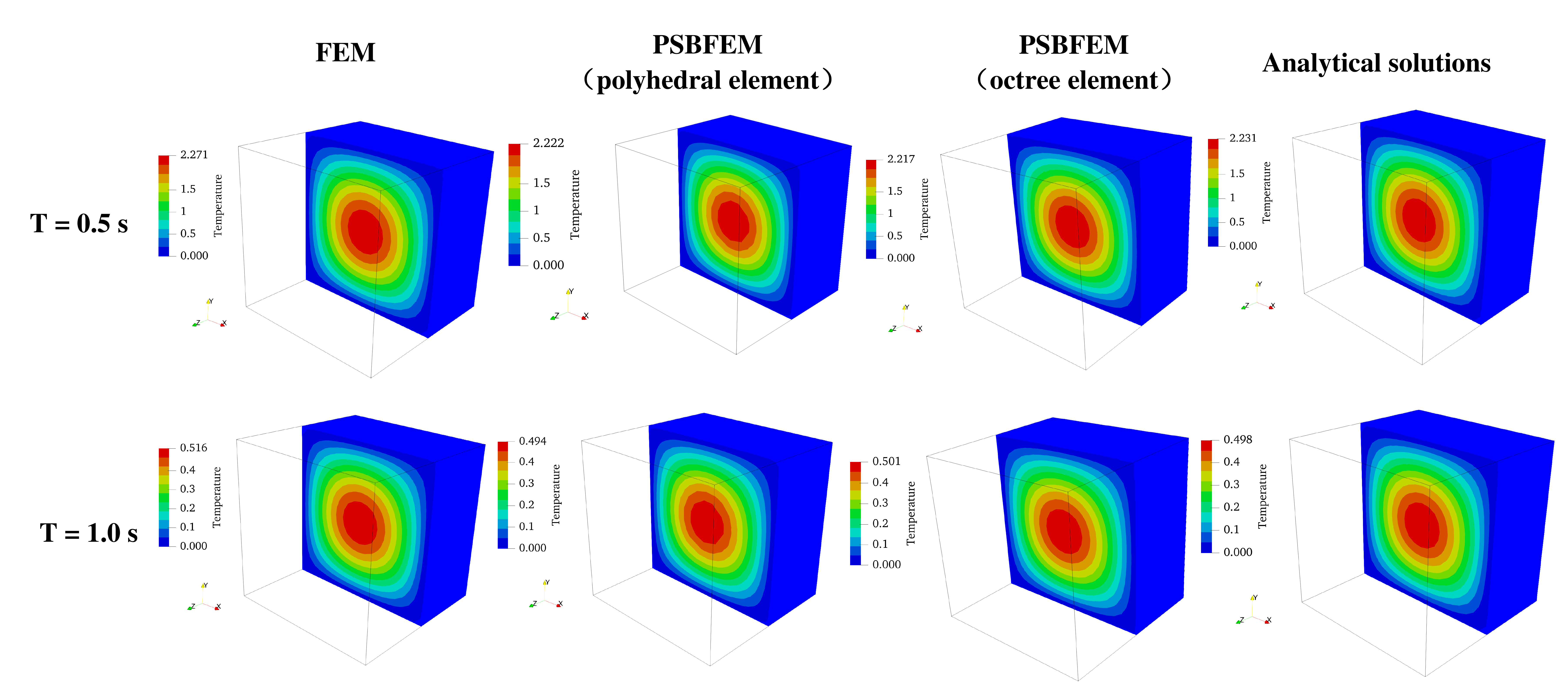}
  \caption{\textcolor{blue}{Distribution of the temperature at different times; (Units:$^{\circ}\mathrm{C}$).}}
  \label{fig:ex03_result}
\end{figure}

\begin{figure}[H]
  \centering
  \includegraphics[width=0.8\textwidth]{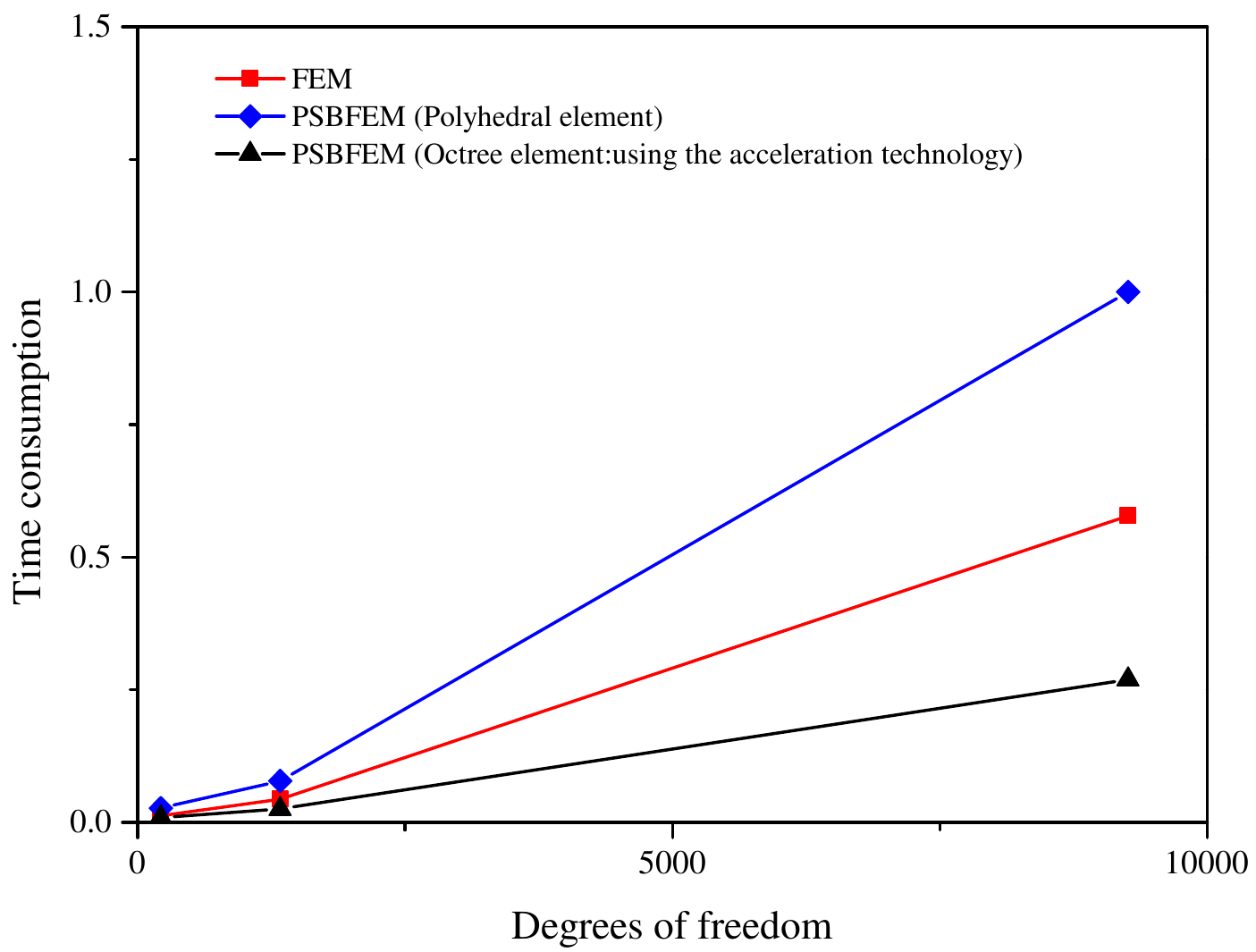}
  \caption{\textcolor{blue}{Time consumption comparison for the transient heat conduction analysis.}}
  \label{fig:ex03_time}
\end{figure}

The comparison of the convergence rates was presented in Fig. \ref{fig:ex03_conv}. It is observed that the convergence rates of the PSBFEM and FEM were identical. At the same element size, the accuracy of the PSBFEM exceeded that of the FEM. The temperature-time history at the monitoring point ($x=\frac{\pi}{2}, y=\frac{\pi}{2}$, and $z=\frac{\pi}{2}$), obtained using the FEM and PSBFEM, was compared in Fig. \ref{fig:ex03_his}, where both methods showed a good correspondence with the analytical solution. Furthermore, Fig. \ref{fig:ex03_result} illustrated the temperature distribution at different time steps for both the PSBFEM and FEM. The results of both methods were in excellent agreement with the analytical solution. \textcolor{blue}{Moreover, for transient problems, the computational time for both methods was shown in Fig. \ref{fig:ex03_time}. Under the same degrees of freedom, the time required by PSBFEM using the polyhedral element was greater than that of FEM, averaging 1.87 times that of FEM. Additionally, with the application of acceleration technology, the computational time for PSBFEM using the octree mesh was reduced to, on average, 58\% of the time required by FEM.} 

\subsection{Complex geometry analysis}
With the rapid advancement of 3D printing, rapid prototyping, and digital analysis, the adoption of the STL format in finite element mesh generation had become increasingly essential. The simplicity and compatibility of the STL format allow efficient representation of complex 3D geometries. In this study, polyhedral meshes and octree-based meshes were utilized to automatically partition STL models, demonstrating the capability of PSBFEM to handle arbitrarily complex polyhedral elements.

\subsubsection{Stanford's bunny using the polyhedral mesh}
\begin{figure}[H]
  \centering
  \includegraphics[width=0.9\textwidth]{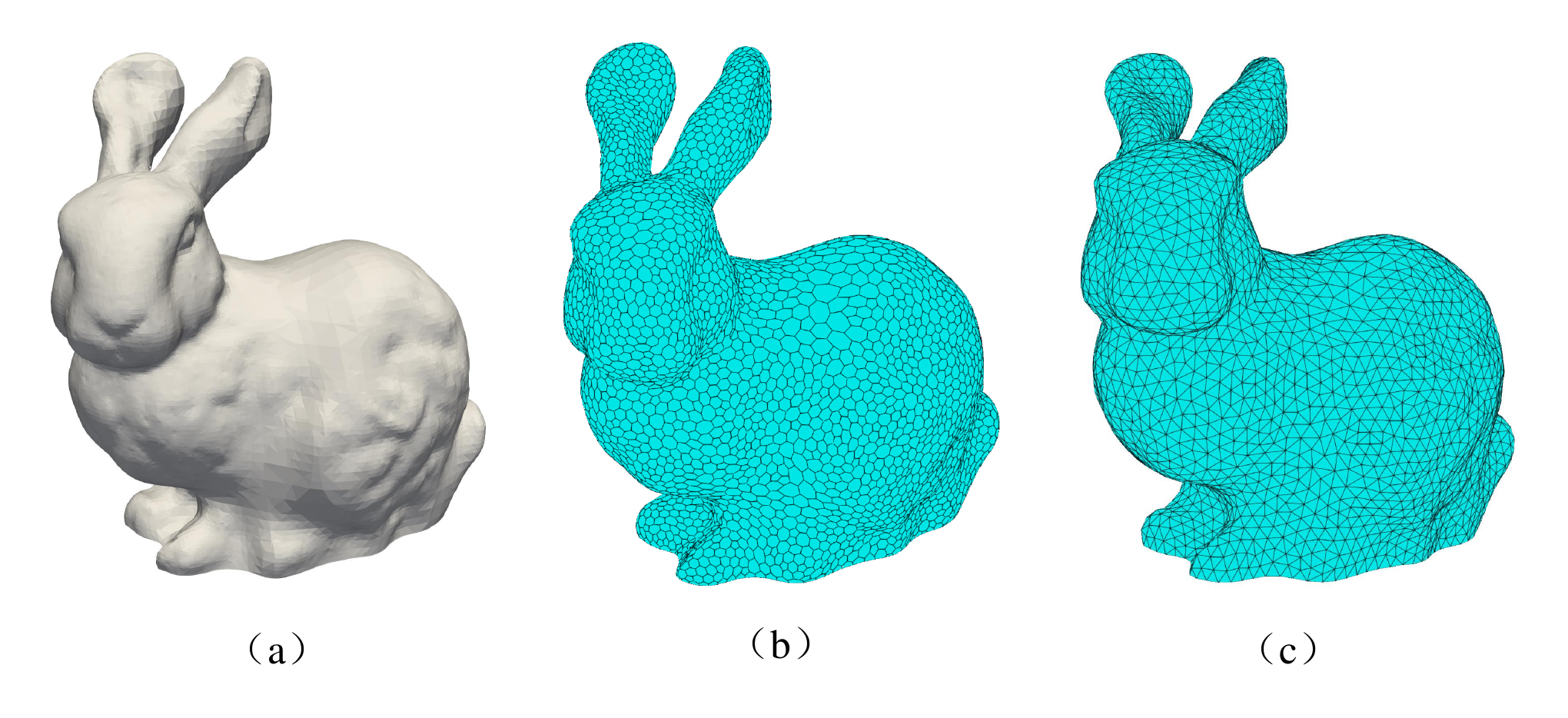}
  \caption{Geometry and mesh model of Stanford's bunny; (a) the geometry model; (b) the polyhedral mesh; (c) the tetrahedral mesh.}
  \label{fig:ex04_bunny_geo_mesh}
\end{figure}

In this example, a Stanford bunny \cite{stanford_bunny} given in STL format was used with a polyhedral mesh \cite{liuAuto2017a}, as shown in Fig. \ref{fig:ex04_bunny_geo_mesh} (b). To further highlight the advantages of the polyhedral mesh, the model was also discretized using an unstructured tetrahedral mesh, as shown in Fig. \ref{fig:ex04_bunny_geo_mesh} (c). The material properties were as follows: thermal conductivity of 52 $\mathrm{W/m/^\circ C}$, specific heat of 434 $\mathrm{J/kg/^\circ C}$, and density of 7800 $\mathrm{kg/m^3}$. A temperature of $T = 1000^\circ \mathrm{C}$ was prescribed on the bottom of the model, while an initial temperature of $T_0 = 20^\circ \mathrm{C}$ was applied throughout the entire model.

\begin{figure}[H]
  \centering
  \includegraphics[width=0.9\textwidth]{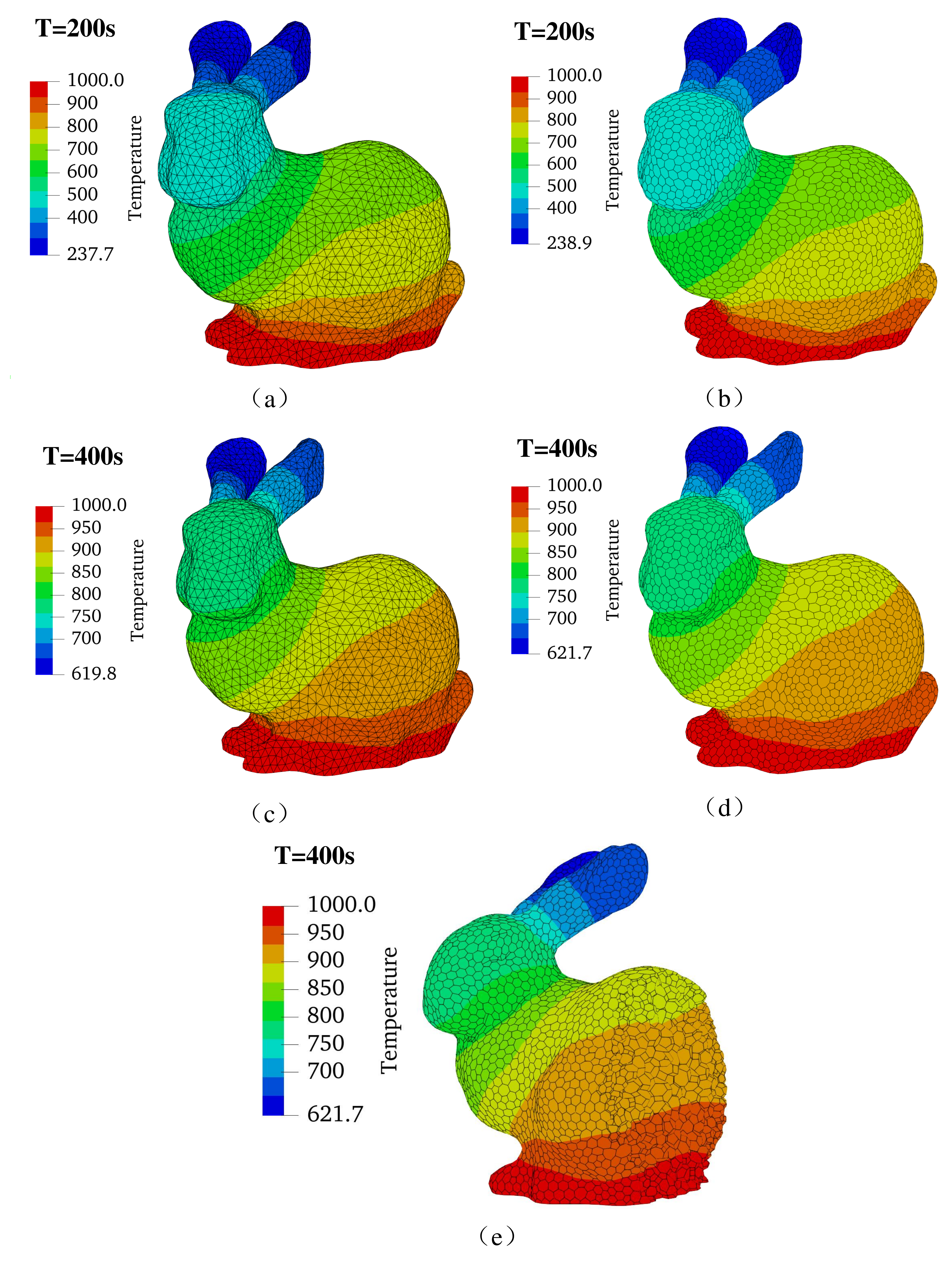}
  \caption{Temperature distribution of Stanford’s bunny using the PSBFEM: (a) temperature contour with tetrahedral elements at 200 s; (b) temperature contour with polyhedral elements at 200 s; (c) temperature contour with tetrahedral elements at 400 s; (d) temperature contour with polyhedral elements at 400 s; (e) localized contour of selected polyhedral elements from (d) at 400 s; (Units:$^{\circ}\mathrm{C}$).}
  \label{fig:ex04_bunnyresult}
\end{figure}

Tab. \ref{tab:Comparison of polyhedron and tetrahedral meshes} presented a comparison between polyhedral and tetrahedral meshes based on several key metrics, including the number of nodes, elements, and surfaces. When using the same element size, the polyhedral mesh consists of 40,714 nodes, whereas the tetrahedral mesh contains only 7,132 nodes. Hence, the polyhedral mesh had a greater number of nodes at the same element size. Additionally, the polyhedral mesh requires only 132.60 seconds of CPU time, compared to the significantly higher 565.20 seconds required by the tetrahedral mesh. This suggested that polyhedral meshes could be a more effective choice for complex simulations, offering reduced computational costs. Furthermore, Fig. \ref{fig:ex04_bunnyresult} illustrated that the temperature distributions for Stanford's bunny using the two different element types exhibited remarkable agreement.

\begin{table}[htbp]
\centering
\caption{Mesh characteristics and relative errors of polyhedral and tetrahedral meshes.}
\begin{tabular}{lcccc}
\toprule
Element type & Nodes & Elements & Surfaces & CPU time (s) \\
\midrule
Polyhedron mesh & 40714 & 7202 & 91879 &  132.60\\
Tetrahedral mesh & 7132 & 32923 & 131692 &  565.20\\
\bottomrule
\end{tabular}
\label{tab:Comparison of polyhedron and tetrahedral meshes}
\end{table}

\subsubsection{Stanford's luck using the hybrid octree mesh}
\begin{figure}[H]
  \centering
  \includegraphics[width=0.9\textwidth]{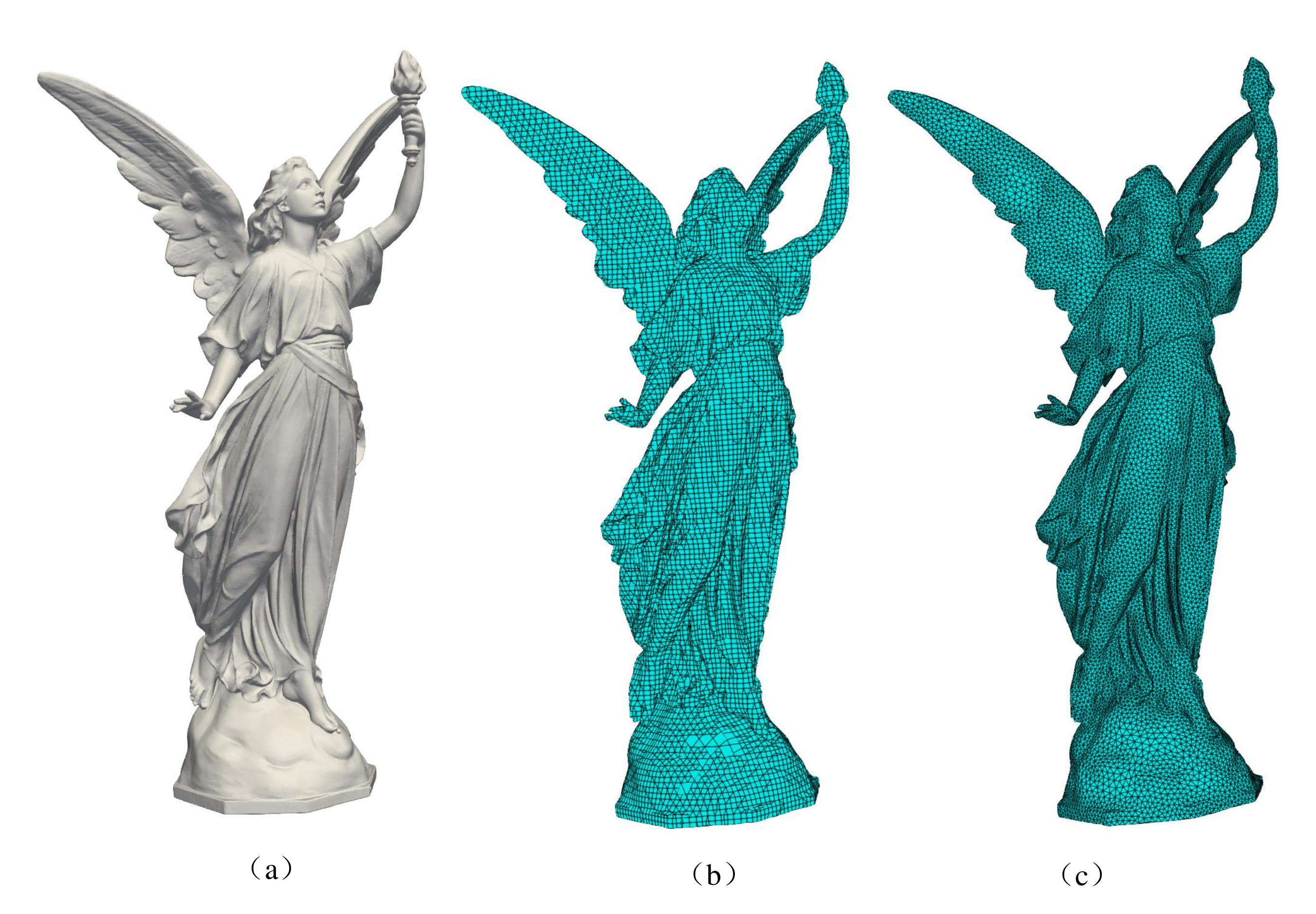}
  \caption{\textcolor{blue}{Geometry and mesh model of Stanford's luck: (a) the geometry model; (b) the hybrid octree mesh; (c) the tetrahedral mesh.}}
  \label{fig:luck_geo_mesh}
\end{figure}

\begin{figure}[H]
  \centering
  \includegraphics[width=0.6\textwidth]{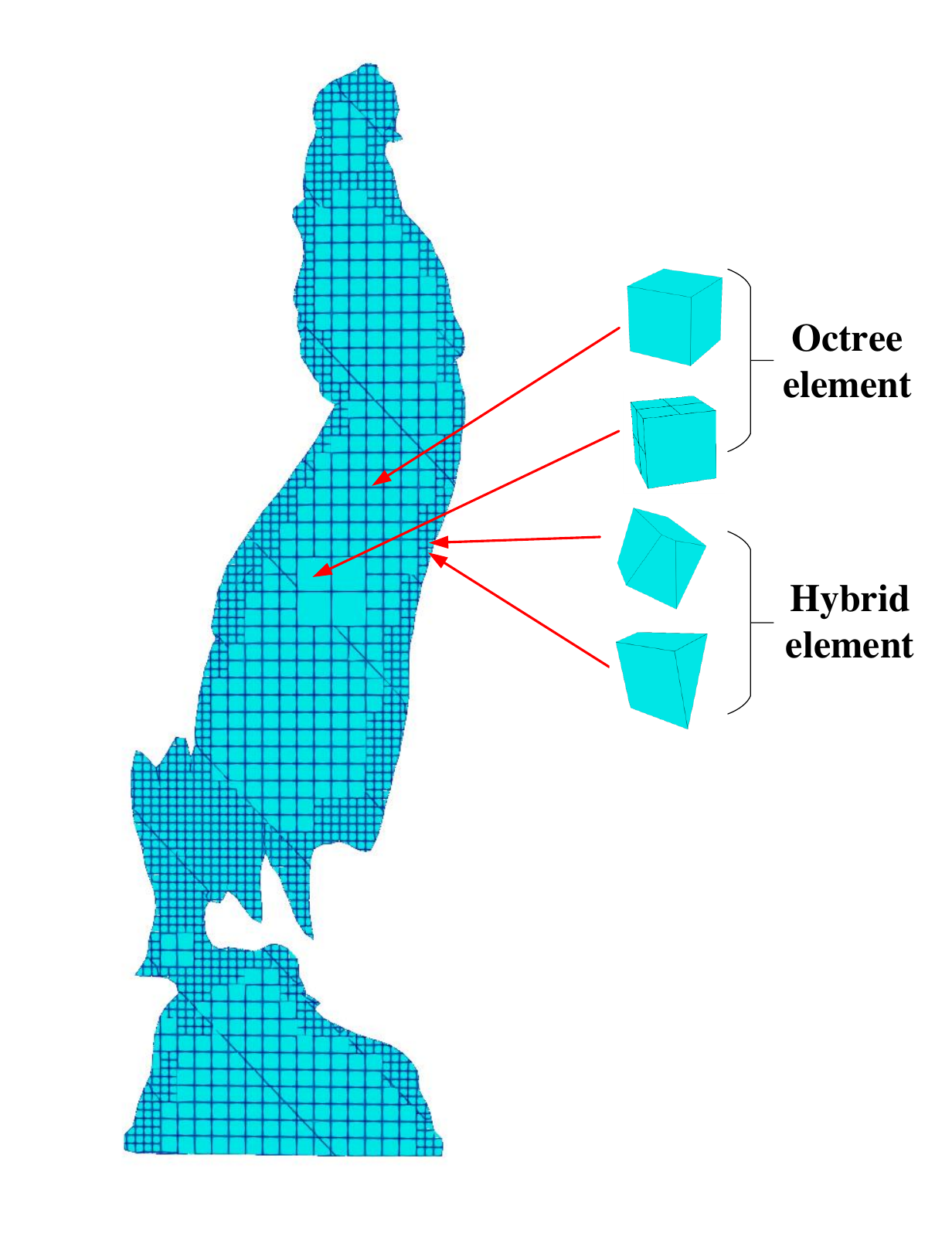}
  \caption{Cross-section view of the hybrid octree mesh for Stanford's luck.}
  \label{fig:ex04_lucksection}
\end{figure}

The hybrid octree mesh was a mesh generation method that combines hexahedral (cubic) elements with other types of elements, such as tetrahedra and polyhedra. For more detailed information, please refer to the relevant literature \cite{songScaledBoundaryFinite2018}. In this example, Stanford's luck \cite{stanford_bunny} was used with the hybrid octree mesh algorithm, as shown in Fig. \ref{fig:luck_geo_mesh} (b). The cross-section of the hybrid octree mesh was illustrated in Fig. \ref{fig:ex04_lucksection}, which showed that the interior of the model was partitioned using a cubic octree mesh, while the outer boundary was formed by generating hybrid elements through a cutting process.

The material properties were defined as: a thermal conductivity of 60 $\mathrm{W/m/^\circ C}$, a specific heat of 97 $\mathrm{J/kg/^\circ C}$, and a density of 2000 $\mathrm{kg/m^3}$. A temperature of $T = 50^\circ \mathrm{C}$ was applied to the bottom of Stanford's luck, while a temperature of $T = 100^\circ \mathrm{C}$ was set on the torch of the model. Tab. \ref{tab:Composition of hybrid octree mesh} presented the composition of the hybrid octree mesh. The mesh was made up of 68.7\% octree elements and 31.3\% hybrid elements. Since the majority of the mesh consisted of high-quality cubic elements, the hybrid octree mesh achieved a high level of accuracy.

\begin{figure}[H]
  \centering
  \includegraphics[width=0.7\textwidth]{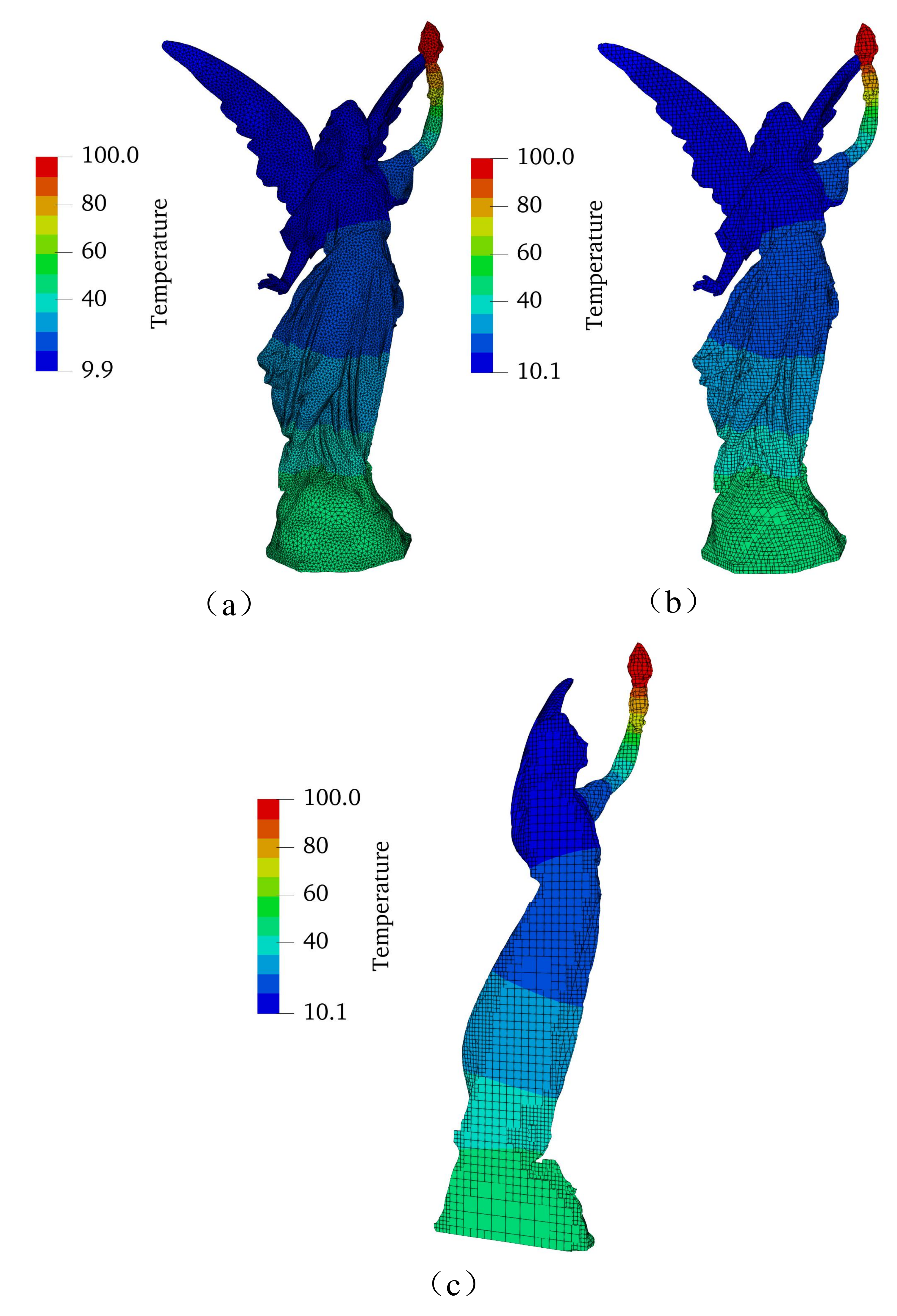}
  \caption{\textcolor{blue}{Temperature distribution of Stanford's luck. (a) temperature contour with tetrahedral elements; (b) temperature contour with hybrid octree elements; (c) localized contour of selected hybrid octree element; (Units:$^{\circ}\mathrm{C}$).}}
  \label{fig:ex04_luckresult}
\end{figure}

Moreover, Tab. \ref{tab:Mesh character hybrid mesh} provided a comparison between the hybrid octree mesh and the tetrahedral mesh. The hybrid octree mesh contained more nodes but fewer elements and surfaces compared to the tetrahedral mesh. The computational cost of the hybrid octree mesh was significantly lower than that of the tetrahedral mesh. \textcolor{blue}{Additionally, through acceleration technology, the computational time was significantly reduced.} Fig. \ref{fig:ex04_luckresult} showed the temperature distribution for Stanford's bunny using both mesh types, demonstrating a high degree of agreement between the results.

\begin{table}[htbp]
\centering
\caption{Composition of hybrid octree mesh for the Stanford's luck.}
\begin{tabular}{lccc}
\toprule
Element type & octree element & hybrid element\\
\midrule
Number of elements &32939	&15020 \\
Proportion &68.7\%	&31.3\% \\
\bottomrule
\end{tabular}
\label{tab:Composition of hybrid octree mesh}
\end{table}


\begin{table}[H]
\centering
\caption{Mesh characteristics, relative errors and time consumption.}
\begin{tabular}{lccccc}
\toprule
Element type &  Nodes & Elements &  Surfaces & CPU time (s) & \textcolor{blue}{CPU time (s)*}\\
\midrule
Hybrid octree mesh & 60579 & 47959 & 293594 & 209.6 & \textcolor{blue}{98.7}\\
Tetrahedral mesh & 48684 & 222210 & 888840 & 859.4 & \textcolor{blue}{—}\\
\bottomrule
\end{tabular}
\label{tab:Mesh character hybrid mesh}
\vspace{-1em}
\begin{flushleft}
\textcolor{blue}{Note: * denotes the use of acceleration technology based on the parent element.}
\end{flushleft}
\end{table}

\section{Conclusions}
\label{sec:8}
This paper presented a PSBFEM framework that demonstrated considerable advantages in solving three-dimensional steady-state and transient heat conduction problems. The key findings of this work are summarized as follows:

(1) By introducing polyhedral elements and utilizing Wachspress shape functions, the proposed method improved mesh efficiency, allowing for accurate solutions with fewer computational elements. This approach reduced preprocessing time, especially for complex geometries where traditional FEM methods would have required extensive re-meshing.

(2) Numerical results demonstrated that the PSBFEM could achieve high accuracy, even with coarse meshes. In particular, the method outperformed traditional FEM in terms of convergence rates, with lower relative errors observed in both steady-state and transient heat conduction problems.

(3) The method’s ability to handle complex geometrical models, such as the Stanford Bunny, further highlighted its versatility. Compared to hexahedral mesh, polyhedral meshes not only provided better accuracy but also significantly reduced computational times, showcasing the potential of this method for large-scale industrial applications.

(4) The hybrid octree mesh structure, in combination with SBFEM, enabled efficient local refinement analysis. This approach offered a balance between accuracy and computational cost, making it suitable for applications requiring high-resolution analysis in localized regions while maintaining global solution stability. \textcolor{blue}{Moreover, by utilizing the octree mesh parent element acceleration technique, the computational efficiency of PSBFEM was significantly enhanced, resulting in a reduction of computation time by 0.53 to 0.58 times compared to the FEM.} 

\textcolor{blue}{In this study, we only considered the use of octree elements without hanging nodes to accelerate the computation. To fully leverage the performance of the octree mesh, future research will consider incorporating octree meshes with hanging nodes as a basic element type.} Moreover, future work could explore extending this framework to other types of physical problems, including thermal stress and multiphysics scenarios. Moreover, implementing adaptive meshing strategies within the PSBFEM framework could further enhance its applicability to real-time engineering simulations.

\section{Acknowledgements}
The National Natural Science Foundation of China (grant NO. 42167046), the Yunnan Province Xing Dian Talent Support Program (grant NO. XDYC-QNRC-2022-0764) and Yunan Funndamental Research Projects (grant NO. 202401CF070043) provided support for this study. 

\appendix
\section{The scaled boundary finite element equation for heat conduction problems}
\label{appendix:A}
\color{blue}
The weighted residual method is employed by multiplying Eq. (\ref{eq:threegov}) with the weighting function $w = w(\xi, \eta, \zeta)$ and integrating over the entire domain. A domain $\Omega$ is considered, where an element on the boundary is scaled within the range $\xi_1 \le \xi \le \xi_2$. Eq. (\ref{eq:threegov}) can be separated into three terms $W_{I}$, $W_{II}$, and $W_{III}$ \cite{song2018scaled,songScaled1999}
\begin{equation}
\begin{aligned}&\underbrace{\int_\Omega w\mathbf{b_1}^\mathrm{T}\mathbf{\tilde{q},_\xi}\mathrm{~d}\Omega}_{W_I}+\underbrace{\int_\Omega w\frac{1}{\xi}(\mathbf{b_2}^\mathrm{T}\mathbf{\tilde{q},_\eta}+\mathbf{b_3}^\mathrm{T}\mathbf{\tilde{q},_\zeta})\mathrm{~d}\Omega}_{W_{II}}\\&
+\underbrace{\mathrm{i}\omega\int_\Omega w\rho c\tilde{T}\mathrm{d}\Omega}_{W_{III}}=0.\end{aligned}\label{eq:ap:threeterms}
\end{equation}

The infinitesimal volume $\mathrm{d}\Omega$ can be expressed as
\begin{equation}
    \mathrm{d}\Omega=\xi^2|\mathbf{J}|\mathrm{d}\xi\mathrm{d}\eta\mathrm{d}\zeta. \label{eq:ap:dv}
\end{equation}

Substituting Eq. (\ref{eq:ap:dv}) into the first and third terms of Eq. (\ref{eq:ap:threeterms}), we obtain the following expressions. The first term, $W_I$, is given by  
\begin{equation}
    W_I=\int_{\xi_1}^{\xi_2}\xi^2\int_{S_\xi}w\mathbf{b}_1^\mathrm{T}\mathbf{\tilde{q}},_\xi|\mathbf{J}|d\eta d\zeta d\xi. \label{eq:ap:wI}
\end{equation}  

Similarly, the third term can be written as  
\begin{equation}
    W_{III}=\int_{\xi_1}^{\xi_2}\mathrm{i}\omega\xi^2\int_{S_\xi} w\rho c\tilde{T}|\mathbf{J}|d\eta d\zeta d\xi.\label{eq:ap:wIII}
\end{equation}

For the second term $W_{II}$, applying Green's theorem to the surface integral over $S^\xi$ converts the integral into a contour integral, resulting in \cite{song2018scaled}
\begin{equation}
    \begin{aligned}
        W_{II} &= \int_{\xi_1}^{\xi_2} \xi \bigg( \oint_{\Gamma^\xi} w \left( q_n^\zeta g^\zeta \,\mathrm{d}\eta + q_n^\eta g^\eta \,\mathrm{d}\zeta \right) \\
        &\quad - \int_{S^\xi} \left( -2w \mathbf{b}_1^\mathrm{T} + w,_\eta \mathbf{b}_2^\mathrm{T} + w,_\zeta \mathbf{b}_3^\mathrm{T} \right) \mathbf{\tilde{q}} |\mathbf{J}| \,\mathrm{d}\eta \,\mathrm{d}\zeta \bigg) \,\mathrm{d}\xi,
    \end{aligned} 
    \label{eq:ap:WII}
\end{equation}
where $q_n^\zeta$ and $q_n^\eta$ are the amplitude of the normal flux on the surfaces. $g^\zeta$ and $g^\eta$ are the magnitude of outward normal vectors to the surfaces.

Substituting Eqs. (\ref{eq:ap:wI}), (\ref{eq:ap:WII}), and (\ref{eq:ap:wIII}) into Eq. (\ref{eq:ap:threeterms}), we obtain  
\begin{equation}
\begin{aligned}
   & \int_{\xi_1}^{\xi_2}\Bigg(\xi^2\int_{S_\xi}w\mathbf{b}_1^\mathrm{T}\mathbf{\tilde{q}}_{,\xi} |\mathbf{J}| \,\mathrm{d}\eta \,\mathrm{d}\zeta 
   + \xi \oint_{\Gamma^\xi} w \left( q_n^\zeta g^\zeta \,\mathrm{d}\eta + q_n^\eta g^\eta \,\mathrm{d}\zeta \right) \\
        &\quad - \xi \int_{S^\xi} \left( -2w \mathbf{b}_1^\mathrm{T} + w_{,\eta} \mathbf{b}_2^\mathrm{T} + w_{,\zeta} \mathbf{b}_3^\mathrm{T} \right) \mathbf{\tilde{q}} |\mathbf{J}| \,\mathrm{d}\eta \,\mathrm{d}\zeta\\
        & + \mathrm{i}\omega\xi^2\int_{S_\xi} w\rho c\tilde{T}|\mathbf{J}|\,\mathrm{d}\eta \,\mathrm{d}\zeta \,\Bigg)\mathrm{d}\xi=0.
    \end{aligned} \label{eq:ap:fourtermNew}
\end{equation}

Eq. (\ref{eq:ap:fourtermNew}) holds when the integrand of the integral over $\xi$ equal to zero \cite{song2018scaled}
\begin{equation}
\begin{aligned}
   & \xi^2\int_{S_\xi}w\mathbf{b}_1^\mathrm{T}\mathbf{\tilde{q}}_{,\xi} |\mathbf{J}| \,\mathrm{d}\eta \,\mathrm{d}\zeta
   + \xi \oint_{\Gamma^\xi} w \left( q_n^\zeta g^\zeta \,\mathrm{d}\eta + q_n^\eta g^\eta \,\mathrm{d}\zeta \right)\\
        &\quad - \xi \int_{S^\xi} \left( -2w \mathbf{b}_1^\mathrm{T} + w_{,\eta} \mathbf{b}_2^\mathrm{T} + w_{,\zeta} \mathbf{b}_3^\mathrm{T} \right) \mathbf{\tilde{q}} |\mathbf{J}| \,\mathrm{d}\eta \,\mathrm{d}\zeta\\
        & + \mathrm{i}\omega\xi^2\int_{S_\xi} w\rho c\tilde{T}|\mathbf{J}|\,\mathrm{d}\eta \,\mathrm{d}\zeta =0.
    \end{aligned} \label{eq:ap:threetermNew}
\end{equation}

The weighting function $w = w(\xi, \eta, \zeta)$ is constructed in the same way as the temperature field in Eq. (\ref{eq:L4})
\begin{equation}
    w(\xi,\eta,\zeta)=\mathbf{N}(\eta,\zeta)w(\xi). \label{eq:ap:wfield}
\end{equation}

By substituting Eq. (\ref{eq:ap:wfield}) into Eq. (\ref{eq:ap:threetermNew}) for an arbitrary \( w(\xi) \), we obtain
\begin{equation}
\begin{aligned}
&\xi^2\int_{S^\xi}\mathbf{B}_1^\mathrm{T}\mathbf{\tilde{q}},_\xi|\mathbf{J}|\mathrm{d}\eta\mathrm{d}\zeta-\xi\int_{S^\xi}(-2\mathbf{B}_1^\mathrm{T}+\mathrm{B}_2^\mathrm{T})\mathbf{\tilde{q}}|\mathbf{J}|\mathrm{d}\eta\mathrm{d}\zeta\\
    &+\mathrm{i}\omega\xi^2\int_{S^\xi}\mathbf{N}(\eta,\zeta)^\mathrm{T}\rho c|\mathbf{J}|\mathrm{d}\eta\mathrm{d}\zeta\\&
    -\xi\oint_{\Gamma^\xi}\mathbf{N}(\eta,\zeta)^\mathrm{T}(q_n^\zeta g^\zeta\mathrm{d}\eta+q_n^\eta g^\eta\mathrm{~d}\zeta)=0.
\end{aligned} \label{eq:ap:threetermnew2}
\end{equation}

Substituting the Eq. (\ref{eq:heat flux}) into Eq. (\ref{eq:ap:threetermnew2}), we obtain
\begin{equation}
\begin{aligned}
    &\xi^2\int_{S^\xi}\mathbf{B}_1^\mathrm{T}\mathbf{k}\left(\mathbf{B}_1\tilde{T},_{\xi\xi}+\frac{1}{\xi}\mathbf{B}_2\tilde{T},_\xi-\frac{1}{\xi^2}\mathbf{B}_2\tilde{T}\right)|\mathbf{J}|\mathrm{d}\eta\mathrm{d}\zeta-\\
    &\xi\int_{S^\xi}(-2\mathbf{B}_1^\mathrm{T}+\mathbf{B}_2^\mathrm{T})\mathbf{k}\left(\mathbf{B}_1\tilde{T},_\xi+\frac{1}{\xi}\mathbf{B}_2\tilde{T}\right)|\mathbf{J}|\mathrm{d}\eta\mathrm{d}\zeta-\\
    &\mathrm{i}\omega\xi^2\int_{S^\xi}\mathbf{N}(\eta,\zeta)^\mathrm{T}\rho c \mathbf{N}(\eta,\zeta)\tilde{T}|\mathbf{J}|\mathrm{d}\eta\mathrm{d}\zeta-\xi\mathbf{F}(\xi)=0.
\end{aligned}
\end{equation}

Introducing the coefficient matrices
\begin{equation}\mathbf{E}_0=\int_{\mathbf{S}_\mathrm{e}}\mathbf{B}_1^\mathrm{T}\mathbf{kB}_1\left|\mathbf{J}\right|\mathrm{d}\eta\mathrm{d}\zeta,\end{equation}
\begin{equation}\mathbf{E}_1=\int_{\mathbf{S}_\mathrm{e}}\mathbf{B}_2^\mathrm{T}\mathbf{kB}_1\left|\mathbf{J}\right|\mathrm{d}\eta\mathrm{d}\zeta,\end{equation}
\begin{equation}\mathbf{E}_2=\int_{\mathrm{S}_\mathrm{e}}\mathbf{B}_2^\mathrm{T}\mathbf{kB}_2\left|\mathbf{J}\right|\mathrm{d}\eta\mathrm{d}\zeta,\end{equation}
\begin{equation}\mathbf{M_0}=\int_{\mathrm{S_e}}\mathbf{N^T}\rho c \mathbf{N}\left|\mathbf{J}\right|\mathrm{d}\eta\mathrm{d}\zeta,\end{equation}
and
\begin{equation}
    \mathbf{F}(\xi)=\int_{\xi_{1}}^{\xi_{2}} [\mathbf{N}(\eta(\xi),\zeta(\xi))]^\mathrm{T} 
    \bigg( q_{n}^{\zeta} g^{\zeta} \frac{\mathrm{d}\eta}{\mathrm{d}\xi} + 
    q_{n}^{\eta} g^{\eta} \frac{\mathrm{d}\zeta}{\mathrm{d}\xi} \bigg) 
    \,\mathrm{d}\xi,
\end{equation}
yields
\begin{equation}
\begin{aligned}
    &\mathbf{E}_0 \xi^2 \tilde{T}(\xi)_{,\xi\xi} 
    + \left( 2\mathbf{E}_0 - \mathbf{E}_1 + \mathbf{E}_1^{T} \right) \xi \tilde{T}(\xi)_{,\xi} \notag \\
    &+ \left( \mathbf{E}_1^\mathrm{T} - \mathbf{E}_2 \right) \tilde{T}(\xi) 
    - \mathbf{M}_0 \mathrm{i}\omega \xi^2 \tilde{T}(\xi)-\xi \mathbf{F}(\xi)=0.  
\end{aligned}
\end{equation}
\color{black}
\bibliographystyle{elsarticle-num} 
\bibliography{cas-refs}





\end{document}